\documentclass[11pt,a4paper]{article}
\usepackage{amsfonts}
\usepackage{amsthm}
\usepackage{amsmath}
\usepackage{amssymb}
\usepackage{mathrsfs}
\usepackage{dsfont}   
\usepackage{latexsym,epsf,graphicx}
\usepackage{caption}
\usepackage[colorlinks=false]{hyperref}
  \hypersetup{
    pdftitle={Robustness of Regularized Pairwise Learning Methods Based on Kernels},
    pdfauthor={\textcopyright DRAFT VERSION.  Andreas Christmann and Ding-Xuan Zhou}
   }
\usepackage[dvips]{color}
\usepackage{url}
\usepackage[sort&compress]{natbib}  
 \bibpunct{(}{)}{;}{a}{,}{,}
\newtheorem{theorem}{Theorem}[section]

\newtheorem{corollary}[theorem]{Corollary}
\newtheorem{definition}[theorem]{Definition}
\newtheorem{lemma}[theorem]{Lemma}
\newtheorem{example}[theorem]{Example}
\newtheorem{assumption}[theorem]{Assumption}

\newcommand{\bc}{\begin{center}}
\newcommand{\ec}{\end{center}}
\newcommand{\bi}{\begin{itemize}}
\newcommand{\ei}{\end{itemize}}
\newcommand{\be}{\begin{equation}}
\newcommand{\ee}{\end{equation}}
\newcommand{\beqna}{\begin{eqnarray*}}
\newcommand{\eeqna}{\end{eqnarray*}}
\newcommand{\bd}{\begin{displaymath}}
\newcommand{\ed}{\end{displaymath}}
\newcommand{\bt}{\begin{tabular}}
\newcommand{\et}{\end{tabular}}

\newcommand{\myem}[1]{\textbf{#1}}

\newcommand{\pen}{\mathrm{pen}}

\newcommand{\rem}[1]{}
\newlength{\fixboxwidth}
\newcommand{\DD}{\mathds{P}}

\newcommand{\N}{\mathds{N}}

\newcommand{\R}{\mathds{R}}
\newcommand{\Rd}{\R^{d}}

\newcommand{\snorm}[1] {\Vert #1 \Vert}

\newcommand{\hnorm}[1] {\Vert #1 \Vert_{\sH}}
\newcommand{\hhnorm}[1] {\Vert #1 \Vert_{\sH}^2}
\newcommand{\inorm}[1] {\Vert #1 \Vert_\infty}

\DeclareMathOperator{\id}{id}

\newcommand{\sL}[2] {\mathcal{L}_{#1}(#2)}
\newcommand{\Lx}[2] {{L}_{#1}(#2)}
\newcommand{\sH}    {H}  
\newcommand{\B}    {\mathcal{B}}   

\def \lb        { \lambda }

\def \e         { \varepsilon }

\def \s         { \sigma }

\def \t         { \tau }
\def \d         { \delta }

\def \x         { \xi }

\newcommand{\fm}     {\Phi}

\def \P           { \mathrm{P} }   
\newcommand{\Pn}  {\P_n}         
\def \Q           { \mathrm{Q} } 

\def \D           { \mathrm{D} }
\newcommand{\PM}  {\mathcal{M}_1}   
\newcommand{\PMXY}  {\mathcal{M}_1(\cX\times \cY)}   

\newcommand{\Ex}{\mathbb{E}}       

\newcommand{\fP}{f_{L,\P,\lb}}
\newcommand{\fLsP}{f_{\Ls,\P,\lb}}

\newcommand{\fT}{f_{L,\D,\lb}}

\def\zeitende{\hfill \quad \hbox{$\vartriangleleft$}}
\def\exampleende{\ifmmode\zeitende\else{\unskip\nobreak\hfil
\penalty50\hskip1em\null\nobreak\hfil\zeitende
\parfillskip=0pt\finalhyphendemerits=0\endgraf}\fi}

\newcommand{\ca}[1]{{\cal #1}}

\newcommand{\Om}{\Omega}

\DeclareMathOperator{\sign}{sign}
\DeclareMathOperator{\Dom}{Dom}

\newcommand{\Int}{\int\limits}

\newcommand{\RPo}[1]{{{\cal R}_{#1,\P}}}

\newcommand{\RPreg}[2]{{{\cal R}_{#1,\P,\lb}^{reg}(#2)}}
\newcommand{\RPs}[2]{{{\cal R}_{#1,\P}'(#2)}}
\newcommand{\RP}[2]{{{\cal R}_{#1,\P}(#2)}}
\newcommand{\RPB}[1]{{{\cal R}_{#1,\P}^{*}}}
\newcommand{\RT}[2]{{{\cal R}_{#1,\D}(#2)}}
\newcommand{\RD}[2]{{{\cal R}_{#1,\D}(#2)}}

\newcommand{\RPr}[2]{{\cal R}_{#1,{\P},\lambda}^{{reg}}(#2)}

\newcommand{\RPL}{RPL }
\newcommand{\RPLM}{RPL method }

\newcommand{\DfiveL}[1] {D_5 L\bigl(X,Y,\tiX,\tiY, {#1}(X), {#1}(\tiX)\bigr)}
\newcommand{\DsixL}[1] {D_6 L\bigl(X,Y,\tiX,\tiY, {#1}(X), {#1}(\tiX)\bigr)}

\newcommand{\DfiveLsmall}[1] {D_5 L\bigl(x,y,\tix,\tiy, {#1}(x), {#1}(\tix)\bigr)}
\newcommand{\DsixLsmall}[1] {D_6 L\bigl(x,y,\tix,\tiy, {#1}(x), {#1}(\tix)\bigr)}

\newcommand{\fLP}{f_{L,\P,\lb}}

\newcommand{\Ls}{L^{\star}}
\newcommand{\fPLs}{f_{\Ls,\P,\lb}}
\newcommand{\fQLs}{f_{\Ls,\Q,\lb}}
\newcommand{\fPepsLs}{f_{\Ls,\P_\e,\lb}}

\newcommand{\fDDnLs}{f_{\Ls,\DD_n,\lb}}

\newcommand{\Lzz} {L_{z,\tiz}}
\newcommand{\Pe} {\P_\e}

\newcommand{\cR}{\mathcal{R}}
\newcommand{\cX}{\mathcal{X}}
\newcommand{\cY}{\mathcal{Y}}
\newcommand{\cZ}{\mathcal{Z}}
\newcommand{\cXY}{{\cX\times\cY}}

\newcommand{\cXYR}{{\cX\times\cY\times\R}}
\newcommand{\tix} {\tilde{x}}
\newcommand{\tiy} {\tilde{y}}
\newcommand{\tiz} {\tilde{z}}
\newcommand{\tiX} {\tilde{X}}
\newcommand{\tiY} {\tilde{Y}}
\newcommand{\tiZ} {\tilde{Z}}
\newcommand{\tit} {\tilde{t}}
\newcommand{\xyxy} {x,y,\tix,\tiy}
\newcommand{\xyxytt} {x,y,\tix,\tiy,t,\tit}
\newcommand{\xyxyfxfx} {x,y,\tix,\tiy,f(x),f(\tix)}

\newcommand{\TT}{{\mathsf{T}}}

\newenvironment{declaration}[1]{\trivlist \item[\hskip \labelsep{\em #1 }]\ignorespaces}{\endtrivlist}
\newenvironment{proofof}[1]{\begin{declaration}{#1.}}{\end{declaration}}

\def \O { \Omega }

\newcommand{\bnum}{\begin{enumerate}}
\newcommand{\enum}{\end{enumerate}}

\newcommand{\IF}{\mathrm{IF}}

\def \lb        { \lambda }

\newcommand{\beq}{\begin{eqnarray}}
\newcommand{\eeq}{\end{eqnarray}}

\newcommand{\XYXY}{X,Y,\tiX,\tiY}

\newcommand{\RLs}[2]{{{\cal R}_{\Ls,#1}(#2)}}
\newcommand{\RLsreg}[2]{{{\cal R}^{reg}_{\Ls,#1,\lb}(#2)}}

\newcommand{\qedr}{\hfill \quad \qed}

\newcommand{\Law}[2]{\mathscr{L}_{#1}(#2)}

\newcommand{\dPro}{{{d}_{\mathrm{Pro}}}} 
\newcommand{\dBL}{{{d}_{\mathrm{BL}}}} 
\def\textregtrademark{\raise1ex\hbox{\scriptsize\textregistered}} 

\newcommand{\NN}{\N}
\newcommand{\RR}{\R}
\newcommand{\E}{{\Ex}}

\numberwithin{equation}{section}

\title{\textbf{On the Robustness of Regularized Pairwise Learning Methods Based on Kernels}$^\dag$\footnotetext{\dag
Corresponding author: A. Christmann, Email: \url{andreas.christmann@uni-bayreuth.de}\newline
~The work by A. Christmann described in this paper is partially supported by a grant of the Deutsche
Forschungsgesellschaft [Project No. CH/291/2-1]. The work by D.-X. Zhou described in this paper is supported partially by a grant from the NSFC/RGC Joint Research Scheme [RGC
Project No. N\_CityU120/14 and NSFC Project No. 11461161006].}}

\author{\textbf{Andreas Christmann}$^1$
and \textbf{Ding-Xuan Zhou}$^2$\\
$^1$ University of Bayreuth, Germany\\
$^2$ City University of Hong Kong, China} 
\date{Date: \today}

\begin{document}

\maketitle

\begin{abstract}
Regularized empirical risk minimization including support vector machines
plays an important role in machine learning theory.
In this paper regularized pairwise learning (RPL) methods based on 
kernels will be investigated. One example is regularized minimization of the error entropy loss which has recently attracted quite some interest from the viewpoint of consistency and learning rates. This paper shows
that such RPL methods have additionally good statistical robustness properties,
if the loss function and the kernel are chosen appropriately.  
We treat two cases of particular interest: (i) a bounded and non-convex loss function and
(ii) an unbounded convex loss function satisfying a certain Lipschitz type condition.
\end{abstract}  
  



\section{Introduction}\label{intro} 

Regularized empirical risk minimization based on kernels have attracted a lot of interest during the last decades in statistical machine learning. 
To fix ideas, let $D=((x_1,y_1), \ldots, (x_n,y_n))$ be a given data set,
where the value $x_i$ denotes the input value and $y_i$ denotes the output value
of the $i$-the data point. 
Let $L$ be a loss function which is typically of the form
$L(x,y,f(x))$, where $f(x)$ denotes the predicted value for $y$, when $x$ is observed, and the
real-valued function $f$ is unknown.
Many regularized learning methods are then defined as minimizers of the optimization problem
\be \label{RERM}
\inf_{f\in\mathcal{F}} \frac{1}{n} \sum_{i=1}^n L(x_i,y_i,f(x_i)) + \pen(\lb,f),
\ee
where the set $\mathcal{F}$ consists of real-valued functions $f$, $\lb>0$ is a regularization
constant, and $\pen(\lb,f)\ge 0$ is some regularization term to ovoid overfitting for the case, that
$\mathcal{F}$ is rich. One example is that $\mathcal{F}$ is a reproducing kernel Hilbert space $H$
and $\pen(\lb,f)=\lb\hhnorm{f}$, see e.g. 
\cite{Vapnik1995,Vapnik1998}, \cite{PoggioGirosi1998}, \cite{Wahba1999}, 
\cite{ScSm2002}, \cite{CuckerZhou2007}, \cite{SC2008} and the references cited there.
 
In recent years there is quite some interest in related learning methods where a 
\emph{pairwise loss function} is used, which yields optimization problems like
\be \label{RPLempirical}
\inf_{f\in H} \frac{1}{n^2} \sum_{i=1}^n \sum_{j=1}^n L(x_i,y_i,x_j,y_j,f(x_i),f(x_j)) + 
\lb \hhnorm{f}
\ee
or asymptotically equivalent versions of it.
In other words, the estimator for $f$ is defined as the minimizer of the sum of 
a $V$-statistic of degree $2$ and the regularizing term $\lb\hhnorm{f}$, see e.g. \cite{Serfling1980}.
An example of this class of learning methods occurs when one is interested in minimizing Renyi's entropy of order 2, see e.g. \cite{HuFanWunZhou2013}, 
\cite{FanHuWuZhou2013}, and \cite{YingZhou2015} for consistency and fast learning rates.
Another example arises from ranking algorithms, see e.g. \cite{Clemencon2008} and \cite{AgarwalNiyogi2009}. Other examples include gradient learning, and metric and similarity learning, see e.g. \cite{MukherjeeZhou2006}, \cite{XinNgJordanRussell2002}, and \cite{CaoGuoYing2015}. 
However, much less theory is currently known for such regularized learning methods 
given by {(\ref{RPLempirical})} based on pairwise loss functions than for the more classical problem {(\ref{RERM})} using standard loss functions. This is true in particular for statistical robustness aspects. Statistical robustness is one important facet of a statistical method,
especially it the data quality is only moderate or unknown, which is often the case in
the so-called big data situation.

The main goal of this paper is to show that such regularized learning methods given by
{(\ref{RPLempirical})} have
nice statistical robustness properties if a combination of a bounded and continuous kernel 
used to define $H$ and a convex, smooth, and separately Lipschitz continuous 
(see Definition \ref{loss:properties}) pairwise loss function is used. 
We also establish a representer theorem for such regularized pairwise learning methods, because we need it for our proofs, but the representer theorem may also be helpful to further research.

The rest of the paper has the following structure.
In Section 2, we define pairwise loss functions, their corresponding risks, derive
some basic properties of pairwise loss functions and their risks, and give some examples.
In Section 3 we define regularized pairwise learning (RPL) methods treated in this paper
and derive results on existence and uniqueness. We will show that shifted loss functions
(defined in {(\ref{shiftedloss})})
are useful to define \RPL methods on the set of \emph{all} 
probability measures without making moment assumptions. This is of course desirable,
because the probability measure chosen by nature to generate the data is completely unknown
in machine learning theory. 
Section 4 contains a representer theorem for \RPL methods, which is our first main result, see
Theorem \ref{thm.representer}. This result is interesting
in its own, but we will also use it as a tool to prove our statistical robustness results in Section 5.
For the case of bounded and not necessarily convex pairwise loss functions
we show that RPL methods have a bounded maxbias if a bounded kernel is used.
For the case of a convex pairwise loss function which is separately Lipschitz continuous in the sense of Definition \ref{loss:properties} and a bounded continuous kernel,
we can formulate the two other main results of this paper: Theorem \ref{thm-boundedderivative}
shows that the \RPL operator has a bounded G\^{a}teaux derivative and hence a bounded influence function, see Corollary \ref{cor-boundedIF}, and Theorem \ref{thm.qualitativerobust} shows that \RPL methods and their empirical bootstrap 
approximations are qualitatively robust, if some \emph{non-stochastic} conditions are satisfied. 
Hence these statistical robustness properties of \RPL methods hold for all probability measures
provided these conditions on the input and output space, on the kernel, and on the loss function
are fulfilled, which can easily be checked by the user.
All proofs are given in the Appendix.

\section{Pairwise Loss Functions and Basic Properties}\label{loss}

If not otherwise mentioned, we will assume the following setup. 

\begin{assumption}\label{assumption-spaces1}
Let $\cX$ be a complete separable metric space and 
$\cY\subset\R$ be closed. 
Let $(X,Y)$ and $(X_i,Y_i)$, $i\in\N$, be independent and identically
distributed pairs of random quantities with values in $\cXY$.
We denote the joint distribution of $(X_i,Y_i)$ by $\P\in\PM(\cXY)$, where $\PMXY$ is the set of all Borel probability measures on the Borel $\s$-algebra $\B_\cXY$. 
\end{assumption}

As usual we will denote the realisations of $(X_i, Y_i)$ by 
$(x_i,y_i)$. For a given data set $D_n=\big( (x_1,y_1),\ldots, (x_n,y_n)\big)\in(\cXY)^n$ we denote the empirical distribution by $\D_n=\frac{1}{n}\sum_{i=1}^n \delta_{(x_i,y_i)}$. Furthermore, we write
$\DD:=\DD_n=\frac{1}{n}\sum_{i=1}^n \delta_{(X_i,Y_i)}$.
Hence $\DD_n(\omega)=\D_n$ for the realisations 
$(X_i(\omega),Y_i(\omega))=(x_i,y_i)$, $i=1,\ldots,n$.

Leading examples on the spaces are of course given by $\cX=\Rd$ and $\cY=\{-1,+1\}$ for binary classification
and $\cX=\Rd$ and $\cY=\R$ for regression, where $d$ is some positive integer.

For $\P\in\PMXY$, we denote the marginal distribution of
$X$ by $\P_\cX$ and the conditional probability of $Y$ given $X=x$
by $\P(y|x)$. The $n$-fold product measure of $\P$ is denoted by
$\P\otimes\ldots\otimes\P$ or simply by $\P^n$.

The classical definition of a loss function in the machine learning literature is a 
measurable function from $\cXYR$ to $[0,\infty)$ and one goal is to minimize
the expected loss plus some regularization term over a hypothesis space, which is often a
reproducing kernel Hilbert space (RKHS), say $H$, defined implicitly by a kernel $k:\cX\times\cX\to\R$,
see e.g. 
\cite{Vapnik1995,Vapnik1998}, \cite{PoggioGirosi1998}, \cite{Wahba1999}, 
\cite{ScSm2002}, \cite{CuckerZhou2007}, \cite{SC2008} and the references cited there.

Here we will consider the case that a regularized risk is to be minimized, where the
loss function for pairwise learning has six instead of three arguments.
I.e., we wish to find a function $f:\cX\to \R$ such that for some non-negative loss function $L$ the value
$L(x,y,\tix,\tiy,f(x),f(\tix))$ is small, if the pair 
$(f(x),f(\tix))$ is a good prediction for the pair $(y,\tiy)$. 
The close connection to $V$-statistics and $U$-statistics (both of degree 2) is obvious, see e.g.
\citet[p.\,172-174]{Serfling1980} and \cite{KoroljukBorovskich1994}.

\begin{definition}\label{loss:u2loss}
Let $(\cX,\ca A)$ be a measurable space and $Y\subset \R$ be closed. 
Then a function
\be \label{def-u2loss}
L:(\cXY)^2\times\R^2 \to [0,\infty) 
\ee
is called a \myem{pairwise loss function}, or simply 
a \myem{pairwise loss}, if it is measurable.
A pairwise loss  $L$ is \myem{represented by $\rho$}, if $\rho:\R\to[0,\infty)$ is
a measurable function and, for all 
$(x,y)\in\cXY$, for all $(\tix,\tiy)\in\cXY$, and for all $t,\tit\in\R$, 
\begin{equation}\label{loss:def-mee}
L(x,y,\tix,\tiy,t,\tit) := \rho\bigl((y-t) - (\tiy-\tit) \bigr).
\end{equation}
\end{definition}

In the following, we will interpret $L(\xyxyfxfx)$
as  the {\em loss} when we predict $(f(x),f(\tix))$ if $(x,\tix)$ 
is observed, but the true outcome is $(y,\tiy)$. 
The smaller the value
$L(\xyxyfxfx)$ is, the better $(f(x),f(\tix))$ predicts $(y,\tiy)$ by means of $L$. 
From this it becomes clear that constant loss functions, such as $L:=0$, 
are rather meaningless for our purposes, since they do not distinguish between good and bad predictions.
Therefore, we will only consider non-constant pairwise loss functions.

Let us now recall from the introduction that 
our major goal is to have a small \emph{average} loss for future unseen observations
$(x,y)$. This leads to the following definition.

\begin{definition}\label{loss:risk-def}
Let $L$ be a pairwise loss function and $\P\in\PM(\cXY)$. 
\bnum
\item The \myem{$L$-risk} for a measurable function $f:\cX\to \R$, i.e. $f\in \sL{0}{\cX}$,  is defined by
\be \label{def-u2risk}
\RP{L}{f}
:=  
\int_{(\cXY)^2} L\bigl(\xyxyfxfx\bigr) \, d\P^2(\xyxy).
\ee
\item The minimal $L$-risk
\be \label{loss:bayesrisk}
  \RPB L := \inf\bigl\{ \RP L f \, ;\,  f:\cX\to \R \mbox{ measurable}   \bigr\}
\ee
is called the \myem{Bayes risk} with respect to $\P$ and $L$. In addition,
a measurable function $f_{L,\P}^{*}:\cX\to \R$ with $\RP L {f_{L,\P}^{*}}= \RPB L$ 
is called a \myem{Bayes decision function}.
\enum
\end{definition}

Note that the function $(x,\tix,y,\tiy) \mapsto L(\xyxyfxfx)$ is measurable by our assumptions, 
and since it is also non-negative,
the above integral always exists, although it is not necessarily finite.

If $\cX$ is a Polish space and $\cY\subset\R$ is closed, then $\cX\times\cY$ is a Polish space, too, such that
we can split up $\P$ into the regular conditional probability $\P(dy|x)$
and the marginal distribution $\P_\cX$, \emph{cf.} \citet[Section 10.2]{Dudley2002}. If we combine this with the
Tonelli-Fubini theorem, we can write $\RP{L}{f}$ as
\be \label{def-u2risk2}
\Int_\cX \Int_\cX \Int_\cY \Int_\cY L\bigl(\xyxyfxfx\bigr) \, \P(dy|x)\,\P(d\tiy|\tix)\,\P_\cX(dx)\,\P_\cX(d\tix).
\ee

For a given sequence $D:=D_n:=((x_1,y_1), \dots, (x_n,y_n))\in (\cX\times \cY)^n$, 
we denote by $\D:=\D_n:=\frac 1 n \sum_{i=1}^{n}\d_{(x_{i},y_{i})}$
the  empirical measure associated to the data set $D$.
The risk of a function $f:\cX\to \R$ 
with respect to $\D$  is  called the \myem{empirical $L$-risk} 
\begin{equation}\label{loss:emp-lrisk}
\RD{L}{f} = \frac{1}{n^2} \sum_{i=1}^{n} \sum_{j=1}^{n} L\bigl(x_i,x_j,y_i, y_j, f(x_i),f(x_j)\bigr)\, .
\end{equation} 
Analogously to {(\ref{loss:emp-lrisk})} one can also define modifications in which one only adds
terms in {(\ref{loss:emp-lrisk})} over all pairs $(i,j)$ with $i\ne j$ or over all pairs $(i,j)$ with
$i\le j$. 
Here we will not investigate these modifications.

We will now introduce some useful properties of pairwise loss functions and their risks 
in much the same manner than these properties are defined for classical loss functions.
The first step is of course measurability.

\begin{lemma}[Measurability of risks]\label{loss:meas-risk}
Let  $L$ be a pairwise loss   and $\ca F\subset \sL 0 \cX$
be a subset that is equipped with a complete and separable metric $d$ and
its corresponding Borel $\s$-algebra.
Assume that the metric $d$ dominates the pointwise convergence, i.e., 
$\lim_{n\to \infty}d(f_n,f) = 0$ implies 
$\lim_{n\to \infty} f_n(x) = f(x)$ for all $x\in \cX$ and 
for all $f, f_n\in \ca F$.
Then the evaluation map 
$\ca F \times \cX \to \R$ defined by $(f,x)\mapsto f(x)$
is measurable, and consequently the map 
$(\xyxy,f)\mapsto L(\xyxyfxfx)$ defined on $(\cXY)^2 \times {\ca F}$ 
and the map 
$(\xyxy,f,f)\mapsto L(\xyxyfxfx)$ 
defined on $(\cXY)^2 \times {\ca F}^2$ are also measurable. 
Finally, given $\P\in\PMXY$, the 
risk functional $\RPo L: \ca F \to [0,\infty]$ is  measurable.
\end{lemma}

Obviously, the metric defined by the supremum norm $\inorm \cdot$
dominates the pointwise convergence for every 
$\ca F\subset C(\cX) \cap \sL \infty \cX$.
It is well-known that the metric of reproducing 
kernel Hilbert spaces (RKHSs) also dominates the 
pointwise convergence.

\begin{definition}\label{loss:properties}
A pairwise loss  $L$ is called  
\begin{enumerate}
\item \myem{(strictly) convex}, \myem{continuous}, or \myem{differentiable}, 
  if 
  $$ L(\xyxy,\,\cdot\,,\,\cdot\,):\R^2\to [0,\infty) $$ 
  is (strictly) convex,
  continuous, or (total) differentiable for all $(\xyxy)\in (\cXY)^2$, respectively.
\item \myem{locally separately Lipschitz continuous}, if  
  for all $b\geq 0$ there exists a constant $c_{b}\geq 0$ such that,
  for all $t,\tit,t',\tit'\in [-b,b]$, we have
  \begin{equation}\label{loss:locallipschitz-loss-def}
  \sup_{\substack{x,\tix\in \cX\\ y,\tiy\in \cY}}  
  \bigl| L(x,y,\tix,\tiy, t,\tit) - L(x,y,\tix,\tiy,t',\tit')\bigr| \,
  \leq \, c_{b} \, \bigl(|t-t'| + |\tit-\tit'|\bigr)\,.
\end{equation}
  Moreover, for $b\geq 0$, the smallest such constant $c_b$ is 
  denoted by $|L|_{b,1}$.
  Furthermore, $L$ is called \myem{separately Lipschitz continuous}\footnote{We mention that \citet{Rio2013}
  used the related term ``separately 1-Lipschitz'' in a different context.}, if
  there exists a constant $|L|_1\in[0,\infty)$ such that,
  for all $t,\tit,t',\tit'\in\R$, the inequality {(\ref{loss:locallipschitz-loss-def})}
  is satisfied, if we replace $c_b$ by $|L|_1$.
\end{enumerate}
\end{definition}

If $L$ is differentiable, we denote by $DL(x,y,\tix,\tiy,t,\tit)$ 
  the (total) derivative of $L(x,y,\tix,\tiy, \,\cdot\,,\cdot\,)$ at $(t,\tit)\in \R^2$.
If $L(x,y,\tix,\tiy,\cdot,\cdot)$ is differentiable with respect to the 5th or 6t
argument, we denote the corresponding partial derivative by
$D_5 L(x,y,\tix,\tiy,\cdot,\cdot)$ and $D_6 L(x,y,\tix,\tiy,\cdot,\cdot)$, respectively.

Let us now consider a few examples of pairwise loss functions.

\begin{example}[Minimum error entropy (MEE) loss]\label{loss:mee}
Fix $h\in(0,\infty)$. Define the pairwise loss 
$L=L_{MEE}$ represented by $\rho_{MEE}$, where 
$\rho_{MEE}(u):=1-\exp\bigl(-u^2/(2h^2)\bigr)$, $u\in\R$,
see e.g. \cite{HuFanWunZhou2013}, \cite{FanHuWuZhou2013}, and
\cite{FengHuangShiYangSuykens2015}.
Some easy calculations show that the first two derivatives $\rho'$ and $\rho''$ are 
continuous and bounded.
However, $\rho_{MEE}$ is not convex and therefore the MEE loss is not convex in the sense of 
Definition \ref{loss:properties}. 
\hfill$\lhd$
\end{example}

\begin{example}[Absolute value type loss]\label{loss:absolutevalue}
Define the pairwise loss  $L$ represented by $\rho_0$, where 
$\rho_0(u)=|u|$, $u\in\R$.
Obviously, $\rho_0$ is Lipschitz continuous, but not differentiable at $u=0$.
\hfill$\lhd$
\end{example}

\begin{figure}
\bc
\includegraphics[width=0.35\textwidth,angle=-90]{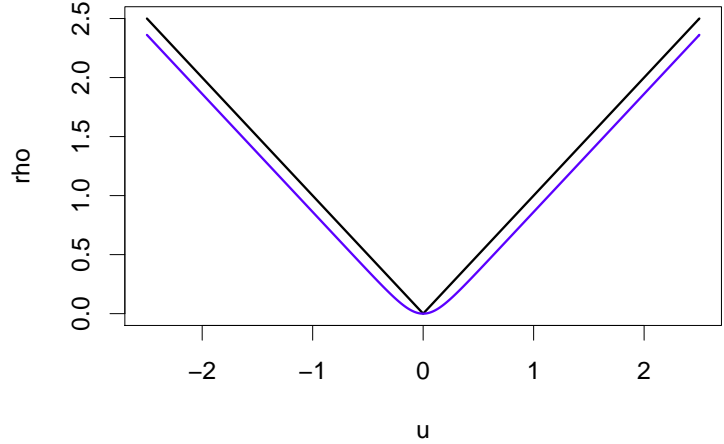} 
 \caption{Comparison of $\rho_0$ (black) and $\rho_a$ with $a=0.01$ (blue).\label{fig:G1}}
\ec
\end{figure}

\begin{example}[Logistic pairwise loss]\label{loss:La}
Fix some $a\in(0,\infty)$, e.g. $a=0.01$ or $a$ equals the rounding precision of the observations.
Denote the cumulative distribution function of the logistic distribution by $\Lambda(r):=1/[1+\exp(-r)]$, $r\in\R$. 
Define the pairwise logistic loss  $L_a$ represented by $\rho_a$, where 
\be \label{rhoa} 
  \rho_a(u):=u - 2 a \log\bigl(2\Lambda(u/a) \bigr), \quad u\in\R.
\ee
As Figure \ref{fig:G1} shows, the pairwise logistic loss 
can be considered for some small tuning value $a$ as a smoothed version of the absolute value type loss.
Some calculations show that $\rho_a$ is Lipschitz continuous with Lipschitz constant $1$ and $\rho_a'$ and $\rho_a''$ are continuous and bounded.
\hfill$\lhd$
\end{example}

\begin{example}[Squared loss]\label{loss:LS}
Define the pairwise loss  $L=L_{LS}$ represented by $\rho_0$, where 
$\rho_{LS}(u)=u^2$, $u\in\R$.
Obviously, $\rho_{LS}$ is only locally Lipschitz continuous and $\rho_{LS}'$ and $\rho_{LS}''$ are
continuous. However, $\rho_{LS}'$ is unbounded.
\hfill$\lhd$
\end{example}

\begin{example}[Ranking loss]\label{loss:ranking} 
Many ranking algorithms can be induced by a pairwise loss of the form 
$L(x,y,\tix,\tiy,t,\tit):=\ell(t-\tit, y-\tiy)$  
with a bivariate function $\ell: \R^2 \to [0, \infty)$. See e.g. \cite{AgarwalNiyogi2009}, 
the hinge ranking loss by 
$$
L(x,y,\tix,\tiy,t,\tit) = \max\{0, v\}
$$
and the least squares ranking loss by
$$
L(x,y,\tix,\tiy,t,\tit) = v^2,
$$
where 
$$
   v:=|y-\tiy| - (t - \tit) \sign(y-\tiy).
$$
For further use we mention that the hinge ranking loss is not differentiable and the least squares ranking loss is not separately Lipschitz continuous in the sense of 
Definition \ref{loss:properties}. 
In contrast, the logistic ranking loss, which we define by
\be 
   L(x,y,\tix,\tiy,t,\tit) = \rho_a(v).
\ee
by using the $\rho_a$ function from {(\ref{rhoa})} for some $a>0$, is a 
separately Lipschitz continuous, differentiable pairwise loss function
with bounded first and second order partial derivatives w.r.t. the last two arguments.
 \hfill$\lhd$
\end{example}

\begin{example}[Similarity loss]\label{loss:distance metric}
Some distance metric or similarity learning algorithms for ${\mathcal X} = \R^d$
and $\cY:=\{-1,+1\}$  can be induced by a pairwise loss of the form 
$$
  L(x,y,\tix,\tiy,t,\tit) = 
  \ell((x -\tix)^T A (x -\tix), r(y, \tiy))
$$ 
or 
$$
  L(x,y,\tix,\tiy,t,\tit) = \ell(x^T A \tix, r(y, \tiy))
$$  
with a positive semidefinite symmetric matrix $A \in \R^{d \times d}$ and $\ell: \R^2 \to [0, \infty)$, $r: \R^2 \to \R$. See e.g. \cite{CaoGuoYing2015}, 
a hinge similarity loss by
$$
  L(x,y,\tix,\tiy,t,\tit) := \max\{0, 1-w\}, 
$$ 
where 
$$ 
   w:=y\tiy x^\TT A \tix.
$$
In the same manner as the logistic pairwise loss function   considered as a smooth alternative to the classical hinge loss function
for binary classification problems, we can define a logistic similarity loss by
\be
  L(x,y,\tix,\tiy,t,\tit) := \ln\bigl( 1+\exp(-w)\bigr). 
\ee
\hfill$\lhd$
\end{example}

It will become clear in Theorem \ref{thm-boundedderivative} and in Theorem \ref{thm.qualitativerobust}, that separate Lipschitz continuity and bounded derivatives are key properties of pairwise loss functions to achieve
a \RPL method  with good robustness properties. 

In the following we will often need that the risk functional is convex to achieve uniqueness of the estimator. This can 
easily be achieved by the following result.

\begin{lemma}[Convexity of risks]\label{loss:convex-loss-lemma}
Let $L$ be a (strictly) convex pairwise loss   and
$\P\in\PMXY$. Then $\RPo L: \sL{0}{\cX}\to [0,\infty]$ is
(strictly) convex.
\end{lemma} 

We also need some additional relationships between a pairwise loss function and its risk.
Such relationships are of course well-known for standard loss functions, see e.g.
\citet{SC2008}.

\begin{lemma}[Lipschitz continuity of risks]\label{loss:lipschitz-loss-lemma}
Let $\P\in\PMXY$ and $L$ be a locally separately Lipschitz continuous pairwise loss. 
Then for all $B\ge 0$ and all $f,g\in L_{\infty}(\P_{\cX})$
with $\inorm{f}\le B$ and $\inorm{g}\le B$, we have
$$
\bigl| \RP{L}{f} - \RP{L}{g} \bigr| 
\le
2 \, |L|_{B,1} \cdot \snorm{f-g}_{\Lx{1}{\P_{\cX}}}\, .
$$ 
Furthermore, the risk functional $\RPo L: L_{\infty}(\P_{\cX})\to [0,\infty)$ 
is well-defined and continuous.
\end{lemma}

In general, we can not expect that the risk of a differentiable loss function is differentiable. 
In this paper we are mainly interested in convex, separately Lipschitz
continuous and differentiable pairwise loss functions, 
for which all partial derivatives (up to order one or two) are continuous and uniformly bounded.
Such loss functions will yield several desirable statistical robustness properties of the learning methods,
as we aim to show. However, we conjecture that
similar results can be shown for certain integrable Nemitski 
losses, see e.g. \citet[Lem.\,2.21]{SC2008} for a result
for standard loss functions from $X\times Y \times \R \to [0,\infty)$.

\begin{lemma}[Differentiability of risks]\label{loss:diff-risk}
Let $\P\in\PMXY$ and $L$ be a differentiable pairwise loss  such that, for all $x,\tix\in\cX$ 
and all $y,\tiy\in\cY$,  the partial derivatives 
$D_i L(x,y,\tix,\tiy,\,\cdot,\,\cdot\,)$, $i\in\{5,6\}$,
are continuous and uniformly bounded by some constant $c_L\in[0,\infty)$.  
Then the risk functional  
$\RPo{L} : \Lx \infty {\P_{\cX}}\to [0,\infty)$ is 
Fr\'echet differentiable and its derivative at 
$f\in \Lx \infty{\P_{\cX}}$ is the 
bounded linear operator 
$\RPs{L}{f} :\Lx \infty{\P_{\cX}}  \to  \R$, where $\RPs{L}{f}{g}$ equals
$$ 
  \Int_{(\cXY)^2} 
      D_5 L(x,y,\tix,\tiy, f(x),f(\tix))g(x) + D_6 L(x,y,\tix,\tiy,f(x),f(\tix)) g(\tix)
      d\P^2(x,y,\tix,\tiy)\, .
$$
\end{lemma}

\begin{example}[Logistic pairwise loss; continuing Example \ref{loss:La}]\label{loss:La2}
For later use, we mention some properties of the pairwise loss  $L_a$, 
where $a\in(0,\infty)$. 
Some tedious but straightforward calculations show that $L_a$ is a convex, 
continuous, differentiable, and separately Lipschitz continuous pairwise loss  with
$|L_a|_1=1$ and 
\beq 
\sup_{(\xyxy) \in (\cXY)^2} | D_i L_a(\xyxy,\,\cdot,\,\cdot\,) | & \le & 1 \label{La:cL1}\\
\sup_{(\xyxy) \in (\cXY)^2} | D_i D_j L_a(\xyxy,\,\cdot,\,\cdot\,) | 
& \le & \frac{1}{2a} \, , \label{La:cL2}
\eeq
where $i,j\in\{5,6\}$. We note that 
$\Lambda'(r)=\Lambda(r)\cdot\bigl(1-\Lambda(r)\bigr) \in(0,1/4)$
and
$\Lambda''(r)=\Lambda(r)\cdot\bigl(1-\Lambda(r)\bigr)\cdot\bigl(1-2\Lambda(r)\bigr)\in(-1/10,+1/10)$, $r\in\R$.
Hence, all partial derivatives of $L_a$ of order up to two w.r.t. the last two arguments are continuous and bounded. 
Fix $a\in(0,\infty)$.
For any fixed values of $(\xyxy)\in(\cXY)^2$ define 
$u:=(y-t)-(\tiy-\tit)=(y-\tilde{y})+(\tit-t)$. 
We have
\beqna
D_5 L_a(\xyxy,t,\tit) & = & -D_6 L_a(\xyxy,t,\tit) ~ = ~ 1-2\Lambda(u/a) \\
D_i D_i L_a(\xyxy,t,\tit) & = & \frac{2}{a} \, \Lambda(u/a) \left( 1- \Lambda(u/a) \right), \qquad i\in\{5,6\}  \\
D_i D_j L_a(\xyxy,t,\tit) & = & -\frac{2}{a} \, \Lambda(u/a) \left( 1- \Lambda(u/a) \right), \quad ~ i,j\in\{5,6\}, i\ne j.
\eeqna
Hence, if we consider $L_a$ as a function of the last two arguments, 
the Hessian matrix is given by
$
\frac{2}{a} \, \Lambda(u/a) \left( 1- \Lambda(u/a) \right) \, \cdot \,
\left[ \begin{array}{rr} +1 & -1 \\ -1 & +1  \end{array} \right]\,,
$
which is obviously positive semi-definite.
Hence, $L_a$ is a convex pairwise loss.
Of course, $L_a$ is a continuous pairwise loss , too.
Furthermore, using the above formulae of the partial derivatives of $L_a$, we obtain that $L_a$ 
is a differentiable pairwise loss  with continuous and bounded partial derivatives
up to order two (w.r.t. the last two arguments of $L_a$).
Therefore, Lemma \ref{loss:diff-risk} yields that
the risk functional 
$\RPo{L} : \Lx \infty {\P_{\cX}}\to [0,\infty)$ is 
Fr\'echet differentiable and its derivative at 
$f\in \Lx \infty{\P_{\cX}}$ is the 
bounded linear operator 
$\RPs{L}{f} :\Lx \infty{\P_{\cX}}  \to  \R$, where  
$\RPs{L}{f}{g}$ is given by
$$
 \Int_{(\cXY)^2} 
 \left(1-2\Lambda\Bigl(\frac{(y-f(x))-(y-f(\tix))}{\e}\Bigr)\right) 
 \cdot \bigl(g(x)-g(\tix)\bigr) \,d\P^2(x,y,\tix,\tiy)\, , 
$$ 
where $g\in \Lx\infty{\P_{\cX}}$. 
Since $1-2\Lambda(r)\in(-1,+1)$ for all $r\in\R$, we immediately
obtain 
$$
| \RPs{L}{f}{g} |
\le 2 \inorm{g}\, , \qquad g\in \Lx\infty{\P_{\cX}}.
$$
Furthermore, some calculations yield that $L_a$ is a separately Lipschitz continuous pairwise loss  with
$|L_a|_1 = 1$ for all $a\in(0,\infty)$.
\hfill$\lhd$
\end{example}

\begin{example}[$L_{LS}$-loss; continuing Example \ref{loss:LS}]\label{loss:LS2}
For later use, we mention that the pairwise loss  $L_{LS}$
is a convex, continuous, and differentiable pairwise loss. 
All partial derivatives of $L_{LS}$ of order up to two w.r.t. the last two arguments are continuous. 
For any fixed values of $(\xyxy)\in(\cXY)^2$ define 
$u:=(y-t)-(\tiy-\tit)=c+\tit-t$. 
We have
\beqna
D_5 L_{LS}(\xyxy,t,\tit) = - D_6 L_{LS}(\xyxy,t,\tit)  & = & -2u \\
D_i D_i L_{LS}(\xyxy,t,\tit) & = & +2,  \qquad i\in\{5,6\} \\
D_i D_j L_{LS}(\xyxy,t,\tit) & = & -2, \qquad i,j\in\{5,6\}, i\ne j.
\eeqna
Hence, if we consider $L_{LS}$ as a function of the last two arguments, 
the Hessian matrix is given by
$
\left[ \begin{array}{rr} +2 & -2 \\ -2 & +2  \end{array} \right]\,,
$
which is obviously positive semi-definite.
Hence, $L_a$ is a convex pairwise loss .
But of course, 
\be \label{LS-unboundedDi}
\sup_{x,\tix\in \cX, ~y,\tiy\in \cY, ~ t,\tit\in\R} | D_i L_{LS}(x,y,\tix,\tiy,t,\tit) | = \infty\,, \quad i\in\{5,6\},
\ee
if $\cY=\R$.
This is in contrast to the separately Lipschitz continuous pairwise loss  $L_a$, as the previous example showed.
\hfill$\lhd$
\end{example}

\section{Regularized Pairwise Learning Methods}\label{existunique}

\begin{definition}\label{def:KMU2}
Let $L$ be a pairwise loss, $H$ be the RKHS of a measurable kernel on $\cX$, and $\lb>0$. For $f\in H$, define the regularized risk by
$\cR_{L,\P,\lb}^{reg}(f) = \RP{L}{f} + \lb \hhnorm{f}$.
A function $\fP \in H$ which satisfies
\begin{equation}\label{infinite:fp-def}
\cR_{L,\P,\lb}^{reg}(\fP) = \inf_{f\in H} \cR_{L,\P,\lb}^{reg}(f)
\end{equation}
is called a \myem{regularized pairwise learning (RLP) method}.
\end{definition}

If $\fP$ exists, we have
\be \label{FPnaiveineq}
\lb\hhnorm{\fP} 
\le 
\cR_{L,\P,\lb}^{reg}(\fP)  
\le
\cR_{L,\P,\lb}^{reg}(0)  
= 
\RP{L}{0}\, ,
\ee
or in other words
\begin{equation}\label{infinite-fp-bound}
\hnorm{\fP} \le \sqrt{\frac{\RP{L}{0}}{\lb}}\, .
\end{equation}

Let us now investigate under which 
assumptions there exists an $\fP\in H$ and when it is unique.  

\begin{assumption}\label{assumption-kernel1} 
Let $k:\cX\times \cX \to \R$ be a \myem{continuous and bounded kernel} 
with reproducing kernel Hilbert space $H$ and define
$\inorm{k}:=\sup_{x\in \cX}\sqrt{k(x,x)}\in(0,\infty)$.
Denote the canonical feature map
by $\Phi(x):=k(\cdot,x)$, $x\in \cX$.
\end{assumption}

It is well-known that if $k$ is a continuous kernel defined on a Polish space, then 
$\Phi$ is continuous, too.
This assumption on the kernel is fulfilled e.g., if $\cX=\Rd$ and $k$ is a Gaussian RBF-kernel
$k(x,x'):=\exp(-\snorm{x-x'}_2^2/\gamma)$, an Abel RBF-kernel
$k(x,x'):=\exp(-\snorm{x-x'}_1/\gamma)$, where $\gamma>0$, or a compactly supported 
kernel, see e.g. \citet{Wu1995} and \citet{Wendland1995}.

The following facts, which we will need later on,  are well-known for any bounded kernel $k$ on $\cX$ with RKHS $H$, all $f\in H$, and 
all $x\in\cX$:
\beq
\langle f, \fm(x) \rangle_H & = & f(x)\,, \label{kernel-prop0}\\
\inorm{f} & \le & \inorm{k} \cdot \hnorm{f}\,, \label{kernel-prop1}\\
\hnorm{\fm(x)} & = & \sqrt{k(x,x)} ~ \le ~ \inorm{k}\, , \label{kernel-prop2}\\
\inorm{\fm(x)} & \le & \inorm{k} \cdot \hnorm{\fm(x)} ~ = ~ \inorm{k}^2\,. \label{kernel-prop3} 
\eeq

First, we will derive conditions for the existence of an \RPL method.

\begin{lemma}\label{weakcong}
Let $r\in(0,\infty)$. 
If $f_0 \in H$ and  if the sequence
$(f_j)_{j=1}^\infty \subset B(f_0, r) :=\{f\in H: \|f - f_0\|_H \leq r\}$, then there exists a subsequence $(f_{j_\ell})_{\ell=1}^\infty$ and $f^{*} \in B(f_0, r)$ such that
 $\|f^{*}\|_H \leq \underline{\lim}_{\ell\to\infty} \|f_{j_\ell}\|_H$ and 
 \begin{equation}\label{converge}
\lim_{\ell\to \infty} f_{j_\ell} (x) = f^{*} (x), \qquad \forall x\in {\mathcal X}. 
\end{equation}
\end{lemma}

\begin{theorem}[Existence]\label{thmexist}
If $L$ is a separately Lipschitz continuous pairwise loss function, $\P\in\PMXY$,
$\cR_{L,\P}(f_0) <\infty$ for some $f_0\in H$, and $H$ be an RKHS with bounded measurable kernel $k$ on $\cX$. 
Then a minimizer $\fP\in H$ exists for any $\lambda >0$.
\end{theorem}

We now address the question of uniqueness.

\begin{lemma}[Uniqueness]\label{infinite:fp-unique}
Let $\P\in\PMXY$, $L$ be a convex pairwise loss with $\RP{L}{f_0}<\infty$ for \emph{some} $f_0\in H$, 
and $H$ be the RKHS of a measurable kernel over $\cX$. 
Then for all $\lb>0$ there exists at most 
one $\fP$.
\end{lemma}

\begin{theorem}[Existence]\label{infinite:fp-exist}
Let $\P\in\PMXY$, $L$ be a convex, locally separately Lipschitz continuous pairwise loss function with 
$\RP{L}{0} < \infty$, and $H$ be the RKHS of a bounded measurable kernel over $\cX$. 
Then, 
for all $\lb>0$, there exists $\fP$.
\end{theorem}

\begin{corollary}[Existence and Uniqueness]\label{infinite:exist-and-unique-mar-based}
Let $\P\in\PMXY$, $L$ be a convex, separately Lipschitz continuous pairwise loss function with 
$\RP{L}{0}<\infty$, and $H$ be the RKHS of a bounded measurable kernel over $\cX$.
Then, for all $\lb>0$,  there exists a uniquely defined $\fP\in H$
and
\be \label{naivefPbound}
\hnorm{\fP} \le \bigl(\RP{L}{0}/\lb\bigr)^{1/2}\, .
\ee
\end{corollary}

Obviously, we would like to get rid of the moment assumption $\RP{L}{0}<\infty$, because otherwise 
we can not define $\fP$ on $\PMXY$ for arbitrary input and output spaces. 
The idea to shift the loss function by an appropriate 
function which is independent of the last two arguments of $L$ is useful in this respect, as was already used e.g. by
\citet{Huber1967} for M-estimators and by \citet{ChristmannVanMessemSteinwart2009} for
support vector machines based on a general loss function $L:\cX\times\cY\times\R\to[0,\infty)$ and on a general kernel.

Let $L$ be a pairwise loss and define the corresponding
\textbf{shifted pairwise loss function} (or simply the shifted version of $L$) by 
\begin{eqnarray} \label{shiftedloss}
   & & \Ls:(\cXY)^2\times\R^2\to\R, \\
   & & \Ls(x,y,\tix,\tiy,t,\tit)  :=  L(x,y,\tix,\tiy,t,\tit)-L(x,y,\tix,\tiy,0,0). \label{shiftedloss2}
\end{eqnarray}
We adopt the definitions of continuity, (locally) separately Lipschitz continuity, and 
differentiability of $\Ls$ from the same definitions for $L$, i.e. these properties are meant to be valid for the last two arguments, when the first four arguments are arbitratily but held fixed.
In the same manner we define the $\Ls$-risk, the regularized $\Ls$-risk, and the \RPLM based on $\Ls$ by
\beq
\RP{\Ls}{f} & := & \Ex_{\P^2} \Ls(X,Y,\tiX,\tiY,f(X),f(\tiX)) \\
\RPreg{\Ls}{f} & := & \RP{\Ls}{f} + \lb\hhnorm{f}\\
\fLsP & := & \arg \inf_{f\in H} \RPreg{\Ls}{f} \,,
\eeq
respectively. 
Of course, shifting the loss function $L$ to $\Ls$ changes the objective function, but
the \emph{minimizers} of $\RPreg{L}{\cdot}$ and $\RPreg{\Ls}{\cdot}$ coincide for those
$\P\in\PM(\cXY)$ for which $\RPreg{L}{\cdot}$
has a minimizer in $H$. I.e.,  we have
\be \label{fLsPequalsfP}
\fLsP = \fP, \mbox{\qquad if~} \fP\in H \mbox{~exists}.
\ee
Furthermore, {(\ref{fLsPequalsfP})} is valid for all empirical distributions $\D$ based on a data set consisting 
of $n$ data points $(x_i,y_i)$, $1 \le i \le n$, because $\fT$ exists and is unique since $\RT{L}{0}<\infty$. 

Let us now show that shifting a pairwise loss function indeed helps to get rid of the moment assumption
$\RP{L}{0}<\infty$ which was essential for the Theorems \ref{infinite:fp-unique}
and \ref{infinite:fp-exist}.
Assume that $L$ is a separately Lipschitz continuous pairwise loss. 
Then we obtain, for all $f \in L_1(\P_\cX)$,
\begin{eqnarray} \label{risk2}
 \RP{\Ls}{f} & =  & \Ex_{\P^2} \bigl( L(X,Y,\tiX,\tiY,f(X),f(\tiX)) - L(X,Y,\tiX,\tiY,0,0) \bigr)  \\
 & \le & \int_{(\cXY)^2} |L(x,y,\tix,\tiy,f(x),f(\tix))-L(x,y,\tix,\tiy,0,0)| \,d\P^2(x,y,\tix,\tiy) \nonumber\\
 & \le & |L|_1 \int_{\cX^2} \bigl( |f(x)-0| + |f(\tix)-0| \bigr) \,d\P_\cX^2(x,\tix) \nonumber \\
 & \le & 2\, |L|_1 \snorm{f}_{L_1(\P_\cX)} ~ < ~  \infty \,, \nonumber
\end{eqnarray}
\emph{without} making the moment condition $\RP{L}{0}<\infty$. 
The assumption  $f \in L_1(\P_\cX)$ can easily be satisfied by choosing a \emph{bounded} kernel $k$, 
because then all $f\in H$ are bounded due to
$\inorm{f} \le \inorm{k}\hnorm{f}$, see {(\ref{kernel-prop1})}.
Therefore, taking {(\ref{fLsPequalsfP})} into account, the use of a shifted loss function just enlarges the set of probability measures where the minimizer of the regularized risk is well-defined. We will make this observation more precise in the remaining part of this 
section.
The following result gives a relationship between $L$ and $\Ls$ in
terms of convexity and Lipschitz continuity.

\begin{lemma} \label{L*prop}
Let $L$ be a pairwise loss. Then the following statements are valid.
\begin{enumerate}
\item If $L$ is (strictly) convex, then $\Ls$ is (strictly) convex.
\item If $L$ is separately Lipschitz continuous, then $\Ls$ is separately Lip\-schitz continuous.
      Furthermore, both Lipschitz constants are equal, i.e., $|L|_1=|\Ls|_1$.
\end{enumerate}
\end{lemma}


\begin{lemma} \label{props}
Let $L$ be a pairwise loss and $\Ls$ its shifted version. Then the following assertions are valid.
\begin{enumerate}
\item $\inf_{f \in \mathcal{L}_0(\cX)} \Ls(x,y,\tix,\tiy,f(x),f(\tix)) \leq 0$.
\item If $L$ is a separately Lipschitz continuous pairwise loss, then for all $f \in H$:
\be \label{riskbound}
  -2 |L|_1 \Ex_{\P_\cX}|f(X)| \le \RP{\Ls}{f} \le 2 |L|_1 \Ex_{\P_\cX}|f(X)|,
\ee
\be \label{riskbound2}  
-2|L|_1 \Ex_{\P_\cX}|f(X)| + \lambda \hhnorm{f}  ~ \le ~  \RPr{\Ls}{f} ~ \le ~ 2 |L|_1 \Ex_{\P_\cX}|f(X)| + \lambda \hhnorm{f}.
\ee
\item $\inf_{f \in H} \RPr{\Ls}{f} \le 0$ and $\inf_{f \in H} \RP{\Ls}{f} \le 0$.
\item Let $L$ be a separately Lipschitz continuous pairwise loss and assume that $\fPLs$ exists.
      Then we have
\begin{eqnarray}
 \lb\hnorm{\fPLs}^2  & \le & - \RP{\Ls}{\fPLs} ~ \le ~ \RP{L}{0}, \nonumber\\
0  & \le & - \RPr{\Ls}{\fPLs} ~ \le ~  \RP{L}{0},\nonumber \\
\label{fbound1} 
\lb\hnorm{\fPLs}^2  & \le &  \min \bigl\{ |L|_1 \Ex_{\P_\cX}|\fPLs(X)|,\RP{L}{0} \bigr\}.
\end{eqnarray}
If the kernel $k$ is additionally bounded, then
\begin{eqnarray} 
\label{fbound3a}  
\inorm{\fPLs} & \le & \lb^{-1} |L|_1 \inorm{k}^2 < \infty,\\
\label{fbound3b}  
|\RP{\Ls}{\fPLs}| & \le & \lb^{-1} |L|_1^2 \inorm{k}^2 < \infty.
\end{eqnarray}
\item If the partial derivatives $D_i L$ and $D_i D_j L$ of $L$ exist for all 
$(x,y,\tix,\tiy) \in (\cXY)^2$ and all $i,j\in\{5,6\}$, then, 
for all $(t,\tit) \in \R^2$,
\begin{eqnarray}  
\label{frecheta1}
D_i \Ls(x,y,\tix,\tiy,t,\tit) & = & D_i L(x,y,\tix,\tiy,t,\tit), \\
\label{frecheta2}
D_i D_j \Ls(x,y,\tix,\tiy,t,\tit) & = & D_i D_j L(x,y,\tix,\tiy,t,\tit).
\end{eqnarray}
\end{enumerate}
\end{lemma}

\noindent The following proposition ensures that the
optimization problem to determine $\fPLs$ is well-posed.

\begin{lemma} \label{proprisk}
Let $L$ be a separately Lipschitz continuous pairwise loss function and $f \in L_1(\P_\cX)$. Then
$\RP{\Ls}{f} \notin \{-\infty, +\infty\}$. Moreover, we have   $\RPr{\Ls}{f} >
-\infty$ for all $f \in L_1(\P_\cX) \cap H$.
\end{lemma}

\begin{lemma}[Convexity of $\Ls$-risks]\label{convrisk2} 
Let $L$ be a (strictly) convex loss. Then
$\RPo{\Ls} : H \to [-\infty,\infty]$ is (strictly) convex and
$\ca{R}_{\Ls,\P,\lb}^{reg} : H \to [-\infty,\infty]$ is strictly convex.
\end{lemma}

\begin{theorem}[Uniqueness of $\fLsP$]\label{infinite:fp-unique2}
Let $L$ be a convex pairwise loss, $H$ be the RKHS of a measurable kernel over $\cX$,
and $\P\in\PMXY$.
Assume that (i) $\RP{\Ls}{f_0}<\infty$ for \emph{some} $f_0\in H$
and $\RP{\Ls}{f}> -\infty$ for \emph{all} $f\in H$
or (ii) $L$ is separately Lipschitz continuous and $f\in L_1(\P_\cX)$ for all $f\in H$. 
Then for all $\lb>0$ there exists at most one decision function $\fLsP$.
\end{theorem}

\begin{theorem}[Existence and Uniqueness of $\fLsP$]\label{infinite:fp-exist2}
Let $L$ be a convex, separately Lipschitz continuous pairwise loss, $H$ be the RKHS of a bounded measurable kernel $k$,
and $\P\in\PMXY$.
Then for all $\lb>0$ there exists a unique decision function $\fLsP$.
\end{theorem}

\section{Representer Theorem for RPL Methods}\label{representer}
In this section we establish a representer theorem for a general probability measure $\P$. This result is  interesting in its own,
but also useful to prove several statistical robustness properties of \RPL methods.
We will often make the following two assumptions to derive 
our representer theorem and robustness results.

\begin{assumption}\label{assumption-loss1} 
Let $L$ be a \textbf{separately Lipschitz-continuous}, \textbf{differentiable} 
pairwise loss function for which all
partial derivatives up to order 2 with respect to the last two arguments are
\textbf{continuous} and \textbf{uniformly bounded} in the sense that there exist
constants $c_{L,1}\in(0,\infty)$ and $c_{L,2}\in(0,\infty)$ with
\beq 
\sup_{x,\tix\in \cX, ~ y,\tiy\in \cY} ~ | D_i L(x,y,\tix,\tiy,\,\cdot,\,\cdot\,) | & \le & c_{L,1}\, , \qquad i\in\{5,6\} \label{loss-assump1}\\
\sup_{x,\tix\in \cX, ~ y,\tiy\in \cY} ~ | D_i D_j L(x,y,\tix,\tiy,\,\cdot,\,\cdot\,) | 
& \le & c_{L,2} \, , \qquad i,j\in\{5,6\} . \label{loss-assump2}
\eeq
Additionally, assume that 
\be \label{Lequalsnull}
  L(x,y,x,y,t,t)\equiv 0, \qquad \forall~ (x,y,t)\in\cX\times\cY\times\R.
\ee
\end{assumption}

Of course, $L_{LS}$ does not satisfy the assumption
{(\ref{loss-assump1})}, whereas e.g. $L_a$ fulfills all conditions in the  Assumption \ref{assumption-loss1}.
The assumption {(\ref{Lequalsnull})} is quite plausible and is satisfied for almost all loss functions
of practical use, e.g. for $L\in\{L_0, L_a, L_{LS}\}$.
If a pairwise loss $L$ is represented by $\rho$, then {(\ref{Lequalsnull})} is satisfied if $\rho(0)=0$.
A ranking loss with $\ell(0,0)=0$ also satisfies {(\ref{Lequalsnull})}.

\begin{assumption}\label{assumption-loss2} 
Let $L$ be a convex pairwise loss function.
\end{assumption}

We will reconsider these assumptions at the end of Section \ref{robust} and it will become clear that these
assumptions on $L$ and $k$ are very plausible to guarantee the existence of 
a \emph{bounded} G\^{a}teaux derivative of the map $\P\mapsto \fLsP$.

As usual we will denote Bochner integrals of an $H$-valued function $g$ with respect to some Borel 
measure $\mu$ by $\int g\,d\mu$, we refer to \citet[p.\, 365ff]{DenkowskiEtAl2003}. If $\mu$ is a probability measure, we denote the Bochner integral occasionally by $\Ex_\mu[g]$. 

\begin{theorem}[\textbf{Representer theorem}\label{thm.representer}]
Let the Assumptions \ref{assumption-spaces1}, \ref{assumption-kernel1}, and \ref{assumption-loss1} be valid.
Then we have, for all $\P\in\PM(\cXY)$ and all $\lb\in(0,\infty)$:
\begin{enumerate}
\item If $\fPLs\in H$ is any fixed minimizer of 
$\min_{f\in H} \big( \RP{\Ls}{f}+\lb\hhnorm{f}\bigr)$, then:
\begin{eqnarray} \label{thm.representer.f1}
 \fPLs  =  - \frac{1}{2\lb} \Ex_{\P^2} 
                       \big[ 
                          h_{5,\P}(\XYXY) \Phi(X) +  h_{6,\P}(\XYXY) \Phi(\tiX)
                       \big],
\end{eqnarray}
where $h_{5,\P}$ and $h_{6,\P}$ denote the partial derivatives
\begin{eqnarray}
   h_{5,\P}(\XYXY)   :=   \DfiveL{\fPLs} ~ \label{thm.representer.f2}\\
   h_{6,\P}(\XYXY)   := \DsixL{\fPLs}. \label{thm.representer.f3}
\end{eqnarray}
\item (Convex Case.) If additionally Assumption \ref{assumption-loss2} is valid (and hence $\fPLs$ uniquely exists by Theorem \ref{infinite:fp-exist2}), then $\fPLs$ has the 
representation {(\ref{thm.representer.f1})} and we have additionally, for all $\Q\in\PM(\cXY)$:
\begin{eqnarray}
 & & \hnorm{\fPLs -\fQLs} \label{thm.representer.f4}\\
 & \le & \frac{1}{\lb} \Big\| 
                      \Ex_{\P^2} 
                       \big[ 
                         h_{5,\P}(\XYXY) \Phi(X)  + h_{6,\P}(\XYXY) \Phi(\tiX)
                       \big]  \nonumber \\
           & & ~~~- \Ex_{\Q^2} 
                       \big[ 
                         h_{5,\P}(\XYXY) \Phi(X)  + h_{6,\P}(\XYXY) \Phi(\tiX)
                       \big]  \Big\|_H \, . \nonumber
\end{eqnarray}
\end{enumerate}
\end{theorem}

\section{Robustness of RPL Methods }\label{robust}
In this section we will show that an \RPL method has several desirable statistical robustness properties, if 
the pairwise loss function $L$ and the kernel $k$ fulfill weak qualitative assumptions. Because these assumptions
are independent of $\P$, these assumptions can really be checked in advance. 
We will start with the case of bounded pairwise loss functions. The case of convex pairwise loss functions 
will be investigated in Section \ref{sec:robustconvex}.

\subsection{Case 1: Non-convex and Bounded Pairwise Loss}\label{sec:robustnonconvex}

The minimizer $f_{\Ls,\P,\lambda}$ typically exists, 
but it is unfortunately in general \emph{not} uniquely defined for \emph{non-convex} pairwise loss functions. 
However we will show
in this subsection, that RPL-methods based on a non-convex
and bounded pairwise loss function often yields a statistically robust
approximations of the  regularized risk. 
More precisely, we will show that the regularized \emph{risk} functional has a small bias
in neighborhoods defined by the norm of total variation, if $L$ is a bounded, but in general \emph{non-convex} 
pairwise loss function. 
This is also valid, if we consider the classical contamination ``neighborhoods'', see
e.g. \citet[p.11]{Huber1981}.
This result will indicate that we can expect a bounded influence function
for the regularized \emph{risk operator} for non-convex pairwise loss functions under appropriate conditions,
provided the influence function exists.

Our most important special case for this section is of course
the minimum entropy loss $L_{MEE}$, see Example \ref{loss:mee}.

In this section, let $L$ be a \emph{bounded} pairwise loss function, 
i.e. we assume 
$$
L(\xyxytt) \in [0,c] \qquad \forall\,(\xyxytt)\in (\cXY)^2 \times\R^2
$$
for some constant  $c\in(0,\infty)$.
Hence the risk $\RP{L}{f}\in[0,c]$ for all $\P\in\PMXY$ and there is no need to consider shifted loss functions. 

The norm of total variation of two probability measures $\P,\Q\in\PM(\cXY)$ is defined by
$$
d_{TV}(\P,\Q):=\sup_{A\in\B(\cXY)}|\P(A)-\Q(A)| = \frac{1}{2} \sup_{h} \left| \int h\,d\P - \int h\,d\Q\right|,
$$
where the supremum is with respect to all $h:\cXY\to\R$ with $\inorm{h}\le 1$. 
It is well-known that $d_{TV}(\P,\Q)\in[0,1]$ for all $\P,\Q\in\PM(\cXY)$.

Define the function 
\begin{eqnarray} \label{generalS}
  & & R^{reg}:\PM(\cXY)\to[0,\infty], \\
  & &  R^{reg}(\P):=\RP{L}{\fLP} = \inf_{f\in H} \RP{L}{f} +  \lb \hhnorm{f}. \nonumber
\end{eqnarray}

Recall that the \emph{maximum bias} of $R^{reg}$ is defined by
\begin{equation}
  b_1(\e;\P):= \sup_{\Q\in N(\e;\P)} |R^{reg}(\Q)-R^{reg}(\P)|, \quad \e\in(0,1),
\end{equation}
where $N(\e;\P)$ denotes an $\e$-neighborhood of $\P$, see \citet[p.11, (4.5)]{Huber1981}.
Common examples are the \emph{total variation neighborhood}
$$  
N_{TV}(\e;\P) :=\{\Q\in\PM(\cXY); d_{TV}(\Q,\P)\le \e\}  
$$ 
and the so-called \emph{contamination ``neighborhood''}
$$ 
N_{con}(\e;\P)=\{\P_\e:=(1-\e)\P+\e\bar{\P}; ~\bar{\P}\in\PM(\cXY)\}.
$$
From a robustness point of view, a statistical method with a bounded maximum bias for sufficiently large positive values
of $\e$ is considered to be robust. If two statistical methods have a bounded maximum bias, the one with the smaller maximum bias is considered to be more robust.

\begin{theorem}[\myem{Bounds for the bias}]\label{boundedmaxbias}
Let $\e\in(0,1)$ and $\P,\Q\in\PM(\cXY)$. Let $L$ be a bounded pairwise loss function satisfying $L\le c\in(0,\infty)$.
Consider the regularised risk functional $R^{reg}$ defined in  {(\ref{generalS})}
\bnum
\item Then  
\begin{equation}
  |R^{reg}(\Q)-R^{reg}(\P)| \le c \, d_{TV}(\Q^2,\P^2) \le 2 \,c \,d_{TV}(\Q,\P)
\end{equation}
and an upper bound for the maximum bias over total variation 
neighborhoods is given by 
$$   b_1(\e;\P)\le 2c\e,$$
uniformly for all $\P$.
\item If $\P_\e=(1-\e)\P+\e\bar{\P}$ for some $\bar{\P}\in\PM(\cXY)$, then  
\begin{equation}
  | R^{reg}(\P_\e)-R^{reg}(\P) | \le 2 \,c \,d_{TV}(\bar{\P},\P) \cdot \e (1+\e)
\end{equation}
and the maximum bias over contamination ``neighborhoods'' satisfies 
$$ b_1(\e;\P)\le 2c\e (1+\e),$$
 uniformly for all $\P$ and $\bar{\P}$.
\enum
\end{theorem}

An obvious consequence of the second part of Theorem \ref{boundedmaxbias} is, that the limit 
\begin{equation}\label{boundedIF0}
\lim_{\e\to 0} \frac{ R^{reg}(\P_\e)-R^{reg}(\P) } {\e}   
\end{equation}
is bounded by $2 \,c \,d_{TV}(\bar{\P},\P) \le 2c$,  
provided the limit exists. If we specialize $\P_\e$ to
$\P_\e=(1-\e)\P+\e\delta_{(x,y)}$ for some $(x,y)\in\cXY$, we obtain immediately from {(\ref{boundedIF0})}, that $S$ has a 
\emph{uniformly bounded} influence function
in sense of \cite{Hampel1968,Hampel1974}, whenever the influence function exists.

Let us now consider an interesting special case of the previous theorem.
Define the discrete probability measures  $\P:=\frac{1}{n}\sum_{i=1}^{n} \delta_{(x_i,y_i)}$ and 
$\P_{n-1}:=\frac{1}{n-1}\sum_{i=1}^{n-1} \delta_{(x_i,y_i)}$ for given data sets
with $n$ and $n-1$ data points $(x_i,y_i)\in\cXY$, respectively, let $\bar{\P}$ be the Dirac measure 
$\delta_{(x_0,y_0)}$ for some point
$(x_0,y_0)\in\cXY$, and let $\e:=\frac{1}{n}$. Then we obtain
$$
\frac{R^{reg}(\P_\e)-R^{reg}(\P)}{\e}
= 
\frac{R^{reg}\bigl((1-\frac{1}{n})\P_{n-1} + \frac{1}{n}\delta_{(x_0,y_0)}\bigr)-R^{reg}(\P_n)}{\frac{1}{n}} \, .
$$
The ratio is the so-called \emph{sensitivity curve} at the point $(x_0,y_0)$, see \citet{Tukey1977} or \citet[p.\,93]{HampelRonchettiRousseeuwStahel1986}, and is usually denoted by
$$SC_n((x_0,y_0);R^{reg},\P_{n-1}).$$
It measures the influence which an additional single data point $(x_0,y_0)$ has on the statistical method $S$,
if the original data set contains $n-1$ data points. The influence function can under appropriate 
assumptions be considered as a finite-sample version of the influence function, see 
\citet[p.\,94]{HampelRonchettiRousseeuwStahel1986}.
A similar version of the sensitivity curve exists, if we
replace one data point from an original data set with $n$ data points.
An immediate consequence of part (ii) in Theorem \ref{boundedmaxbias} is, that
the sensitivity curve $SC_n((x_0,y_0);R^{reg},\P_{n-1})$ is \emph{uniformly bounded} by
$2c(1+\frac{1}{n})$ for all data sets and any additional data point $(x_0,y_0)$, no matter where 
$(x_0,y_0)$ is located in $\cXY$. If we are interested in the slightly more general problem 
how to obtain an upper bound for the influence of $\ell$ additional data
points, we just define $\e:=\frac{\ell}{n}$, $\ell\in\{1,\ldots,\lfloor\frac{n}{2}\rfloor\}$,
and use again part (ii) in Theorem \ref{boundedmaxbias} to obtain a uniform upper bound.

\begin{example}\label{examplemaxbias}
Theorem \ref{boundedmaxbias} is applicable for the non-convex minimum entropy loss $L_{MEE}$ represented by $\rho_{MEE}(u)\in[0,1)$, $u\in\R$, where $h\in(0,\infty)$, see Example \ref{loss:mee}. 
A division by $\e$ shows that the absolute value of these difference quotients are  bounded by $2$ or $2(1+\e)$, respectively, which is an immediate consequence of 
Theorem \ref{boundedmaxbias}.
If we additionally assume for the case of contamination ``neighborhoods'' in 
Theorem \ref{boundedmaxbias}\emph{(ii)}, that the limit
$\lim_{\e \searrow 0} \frac{R^{reg}(\Pe) - R^{reg}(\P)}{\e}$
exists, where $\Pe:=(1-\e)\P+\e\delta_z$ with $\delta_{(x,y)}$ being the Dirac measure in $(x,y)\in\cXY$, 
then this limit equals the influence function of $R^{reg}$ at $\P$
and its absolute value is then bounded by $\lim_{\e\searrow 0}\frac{2\e(1+\e)}{\e}=2$.
This is of course desirable from a robustness point of view.
\end{example}

Summarizing, we showed in this subsection that the regularized \emph{risk} functional based on 
bounded pairwise loss functions $L$ has some desirable robustness properties even if $L$ is \emph{non-convex} and bounded.
An example is the minimum error entropy loss.

\subsection{Case 2: Convex Pairwise Loss} \label{sec:robustconvex}



An immediate consequence of the second part of our representer theorem, see {(\ref{thm.representer.f4})}, is the inequality
\be 
\hnorm{\fPLs -\fQLs} \le \frac{4}{\lb} c_{L,1} \inorm{k}^2 < \infty,
\ee
which is valid for all $\P,\Q\in\PM(\cXY)$ and all $\lb\in(0,\infty)$ if the Assumptions \ref{assumption-spaces1}, \ref{assumption-kernel1}, \ref{assumption-loss1}, and \ref{assumption-loss2} are valid.

The goal of this subsection however is to show that the \textbf{RPL operator}
\be \label{f.RPLoperator}
  S: \PM(\cXY)\to H, \qquad S(\P)=\fLsP
\ee
has two additional desirable robustness properties, if weak conditions on $\cX$, $\cY$, $L$, and $k$ are satisfied: 
\bnum
\item Theorem \ref{thm-boundedderivative} will show that $S$ has a \emph{bounded} G\^{a}teaux derivative for any probability measure $\P$ and hence a \emph{bounded} influence function in the sense
of \cite{Hampel1968,Hampel1974}, see also \cite{HampelRonchettiRousseeuwStahel1986}.
\item Theorem \ref{thm.qualitativerobust} will show that the sequence of RPL estimators
      $(\fDDnLs)_{n\in\N}$ are qualitatively robust, which is a kind of equicontinuity described later in more detail. If additionally $\cXY$ is a compact metric space, then even the empirical bootstrap approximations are qualitatively robust.
\enum

Please note, that the following results of this subsection are all formulated for $\fLsP$ and not for $\fP$, because the latter is in general not well-defined for \emph{all} $\P\in\PMXY$, as was explained in Section \ref{existunique}. Please recall the obvious equalities $D_i L(\xyxytt)=D_i \Ls(\xyxytt)$ and
$D_i D_j L(\xyxytt)=D_i D_j \Ls(\xyxytt)$ for $i,j\in\{5,6\}$.

\begin{theorem}[\myem{Bounded  G\^{a}teaux derivative}]\label{thm-boundedderivative}
Let the Assumptions \ref{assumption-spaces1}, \ref{assumption-kernel1}, \ref{assumption-loss1}, and \ref{assumption-loss2} be satisfied.
Denote the shifted version of $L$ by $\Ls$. 
Then, for all Borel probability measures $\P, \Q\in\PMXY$, 
the operator 
$$
S: \PMXY \to H, \qquad S(\P) := \fLsP
$$
has a bounded G\^{a}teaux derivative $S'_G(\P)$ at $\P$ and it holds
\be \label{thm-Gateaux}
S'_G(\P)(\Q) = -M(\P)^{-1} T(\Q;\P)\,.
\ee
Here 
\beqna
T(\Q;\P) & = & -2 \, \Ex_{\P\otimes\P} \Bigl[ D_5 L\bigl(X,Y,\tiX,\tiY,\fLsP(X),\fLsP(\tiX)\bigr) \Phi(X) \\
         & & ~~~~~~~~~~~~      + D_6  L\bigl(X,Y,\tiX,\tiY,\fLsP(X),\fLsP(\tiX)\bigr) \Phi(\tiX) \Bigr] \\
         & & +~~ \Ex_{\P\otimes\Q} \Bigl[ D_5 L\bigl(X,Y,\tiX,\tiY,\fLsP(X),\fLsP(\tiX)\bigr) \Phi(X) \\
         & & ~~~~~~~~~~~~      + D_6  L\bigl(X,Y,\tiX,\tiY,\fLsP(X),\fLsP(\tiX)\bigr) \Phi(\tiX) \Bigr] \\
         & & +~~ \Ex_{\Q\otimes\P} \Bigl[ D_5 L\bigl(X,Y,\tiX,\tiY,\fLsP(X),\fLsP(\tiX)\bigr) \Phi(X) \\
         & & ~~~~~~~~~~~~      + D_6  L\bigl(X,Y,\tiX,\tiY,\fLsP(X),\fLsP(\tiX)\bigr) \Phi(\tiX) \Bigr] 
\eeqna
equals the gradient of the regularized risk and  
\beqna
M(\P) & = & ~~\, 2 \lb \id_H \\
      & & + \,  \Ex_{\P\otimes\P}  \Bigl[ D_5 D_5 L\bigl(X,Y,\tiX,\tiY,\fLsP(X),\fLsP(\tiX)\bigr) \Phi(X) \otimes \Phi(X) \\
      & & + ~~~~~~~~\,  D_6 D_5 L\bigl(X,Y,\tiX,\tiY,\fLsP(X),\fLsP(\tiX)\bigr) \Phi(X) \otimes \Phi(\tiX) \\
      & & + ~~~~~~~~\,  D_5 D_6 L\bigl(X,Y,\tiX,\tiY,\fLsP(X),\fLsP(\tiX)\bigr) \Phi(\tiX) \otimes \Phi(X) \\
      & & + ~~~~~~~~\,  D_6 D_6 L\bigl(X,Y,\tiX,\tiY,\fLsP(X),\fLsP(\tiX)\bigr) \Phi(\tiX) \otimes \Phi(\tiX) \Bigr]
\eeqna
equals the Hessian of the regularized risk.
\end{theorem}

Please note, that the operator $M(\P)$ and the first integral of $T(\Q;\P)$ only depend on $\P$.
Only the second and the third integral in the formula of $T(\Q;\P)$ depend on $\Q$ and describe how the  $f_{L,\bullet,\lb}$ 
changes, if the probability measure equals the mixture  $(1-\e)\P+\e\Q$
instead of $\fP$. Of course, we have $S'_G(\P)(\P)=0\in H$.

The influence function is an important approach in robust statistics and was 
proposed by \citet{Hampel1968,Hampel1974}; we refer also to the classical textbook by 
\cite{HampelRonchettiRousseeuwStahel1986}. The influence function is related to 
G\^{a}teaux differentiation of the operator $S$ in direction of the Dirac measure 
$\Q:=\delta_{(x_0,y_0)}$, where $(x_0,y_0)\in \cXY$, i.e.
$$
\IF((x_0,y_0); S,\P)=S'_G(\P)(\delta_{(x_0,y_0)}).
$$


The influence function has the interpretation that it measures the influence of an (infinitesimal) small amout of contamination
of the original measure $\P$ in the direction of a Dirac measure located in the point
$(x_0,y_0)$ on the theoretical quantity $S(\P)$ of interest.  Hence, it is desirable that a 
statistical method has a \emph{bounded} influence function. If different methods have a bounded
influence function, the one with the lower bound is considered to be more robust within this approach.

\begin{corollary}[\myem{Bounded influence function}]\label{cor-boundedIF}
Let the Assumptions \ref{assumption-spaces1}, \ref{assumption-kernel1}, 
\ref{assumption-loss1}, and  \ref{assumption-loss2} be satisfied. 
Denote the shifted version of $L$ by $\Ls$. 
Then, for all $\P\in\PMXY$, for all $(x_0,y_0)\in \cX\times \cY$, and
for all $\lb\in(0,\infty)$,
the influence function of
$S: \PMXY \to H$ defined by $S(\P):= \fLsP$ is bounded. It holds
\be \label{robust-formulaIF}
{\rm{IF}}( (x_0,y_0); S, \P) = -M(\P)^{-1} T(\d_{(x_0,y_0)};\P),
\ee
where $T(\d_{(x_0,y_0)};\P)$ and $M(\P)$ are given by Theorem \ref{thm-boundedderivative}
and $T(\d_{(x_0,y_0)};\P)$ simplifies to
\beqna
& & - 2 \, \Ex_{\P\otimes\P} \Bigl[ D_5 L\bigl(X,Y,\tiX,\tiY,\fLsP(X),\fLsP(\tiX)\bigr) \Phi(X) \\
& & ~~~~~~~~~~~~      + D_6  L\bigl(X,Y,\tiX,\tiY,\fLsP(X),\fLsP(\tiX)\bigr) \Phi(\tiX) \Bigr] \\
& & +~~ \Ex_{\P} \Bigl[ D_5 L\bigl(X,Y,x_0,y_0,\fLsP(X),\fLsP(x_0)\bigr) \Phi(X) \\
& & ~~~~~~~~~          + D_6  L\bigl(X,Y,x_0,y_0,\fLsP(X),\fLsP(x_0)\bigr) \Phi(x_0)  \\
& & ~~~~~~~~~          + D_5 L\bigl(x_0,y_0,X,Y,\fLsP(x_0),\fLsP(X)\bigr) \Phi(x_0) \\
& & ~~~~~~~~~          + D_6  L\bigl(x_0,y_0,X,Y,\fLsP(x_0),\fLsP(X)\bigr) \Phi(X) \Bigr]\,.
\eeqna
\end{corollary}

We mention that the pairwise loss $L_a$ fulfills the Assumptions \ref{assumption-loss1} and \ref{assumption-loss2} with
$c_{L_a, 1}=1$ and $c_{L_a,2}=\frac{1}{2a}$ for any $a\in(0,\infty)$, see
{(\ref{La:cL1})} and {(\ref{La:cL2})}. 
Hence, Theorem \ref{thm-boundedderivative} and Corollary \ref{cor-boundedIF}
are applicable for $L_a$, if used in combination with a bounded and continuous kernel, 
e.g. a Gaussian RBF kernel.

Now let us reconsider the assumptions on $L$ and $k$ we made to establish 
Theorem \ref{thm-boundedderivative}. Due to 
$S'_G(\P)(\Q) = -M(\P)^{-1} T(\Q;\P)$ and the specific form of $T(\Q;\P)$ and $M(\P)$, we see that the boundedness of the
G\^{a}teaux derivative stems from the fact that $L$ is separately Lipschitz continuous \emph{and} $k$ is bounded. One of these
properties will in general not be enough to guarantee the boundedness of the G\^{a}teaux derivative in unbounded spaces $\cX$ and $\cY$.
Let us give one simple example. If $\cY$ is unbounded, e.g. $\cY=\R$, we do not expect a bounded influence function for $\fLsP$, 
if the squared loss $L_{LS}$ is used, because the supremum of the absolute values of the partial derivatives are unbounded in this case, 
as follows from {(\ref{LS-unboundedDi})}. Please note that this is no contraction to Theorem \ref{thm-boundedderivative}, because $L_{LS}$ is clearly not separately Lipschitz continuous
and $f_{L_{LS}^\star,\P,\lb}$ is in general not even defined on the set of \emph{all} Borel probability measures
$\PM(\cXY)$, if $\cY$ is unbounded.

In this sense, Theorem \ref{thm-boundedderivative} and its corollary are in good agreement with results obtained by 
\citet{ChristmannSteinwart2004a,ChristmannSteinwart2007a} for the case of support vector machines based on a general loss function
and on a general kernel.


Besides the maximum bias over neighbourhoods and a bounded influence function, qualitative robustness is another key notion in robust statistics. Qualitative robustness was proposed by \citet{Hampel1968,Hampel1971} and generalized to more abstract spaces by \citet{Cuevas1988}. 
Define
$$
  \DD_n:=\frac{1}{n} \sum_{i=1}^n \delta_{(X_i,Y_i)}.
$$
In this subsection we will show that the sequence of estimators 
$$ 
  (S_n)_{n\in\N}, \quad S_n:=f_{\Ls,\DD_n,\lb}
$$ 
is qualitatively robust for all probability measures and any fixed regularization parameter $\lb\in(0,\infty)$.
We will also give an analogous qualitative robustness result for the empirical bootstrap approximations.

According to \cite{Hampel1968} and \citet{Cuevas1988}
a sequence of estimators $(S_n)_{n\in\N}$ is called 
\emph{qualitatively robust} at a probability measure $\P$ if and only if 
\begin{eqnarray}  \label{def.qualitativerobustness}
 \forall\, \varepsilon>0 ~~ \exists\, \delta>0~  : ~~~
 \Bigl[ d_*(\Q,\P)\,<\,\delta \quad \Longrightarrow \quad
  d_*\big(\Law{\Q}{S_{n}}, \Law{\P}{S_{n}}\big)
  \,<\,\varepsilon
    \;\;\forall\,n\in\N \Bigr]\;.
\end{eqnarray}
Here $\Law{\Q}{S_{n}}$ and $\Law{\P}{S_{n}}$ denote the image measures of 
$\Q$ and $\P$ by $S_n$, if all pairs $(X_i,Y_i)$ are independent and identically distributed 
with $(X_i,Y_i) \sim \P$ or $(X_i,Y_i)\sim \Q$, respectively. Another common
notation for $\Law{\P}{S_{n}}$ is $S_n(\P^n)$.
Originally, \citet{Hampel1971} used for $d_{*}$ the Prohorov metric $\dPro$, but one can also use 
the bounded Lipschitz metric $\dBL$ defined by
$$ 
  \dBL(\P,\Q):=\sup \left\{ \left| \int g \,d\P - \int g \,d\Q \right| ; \| g\|_{\mathrm{BL}} \le 1 \right\} \, ,
  ~~ \P,\Q\in\PM(\cXY)
$$
in separable metric spaces, where 
$\|g\|_{L}:= \sup_{x_1 \ne x_2} |g(x_1)-g(x_2)|/d(x_1,x_2)$ and
$\|g\|_{\mathrm{BL}}:=\|g\|_{L} + \inorm{g}$,
see \citet[Chapter 11.2]{Dudley2002}.
The reason for this is that, for any \emph{separable} metric space -- and in our case $\cXY$ is separable --, both $\dPro$ and $\dBL$ metrize the weak convergence for sequences of probability measures, i.e. 
\be \label{WeakConvergenceAndDBL} 
  \P_n \rightsquigarrow \P  \quad \Longleftrightarrow \quad \dPro(\Pn,\P)\to 0 \quad \Longleftrightarrow \quad \dBL(\P_n,\P)\to 0,
\ee
we refer to \citet[Thm. 11.3.3, p. 395]{Dudley2002} for details. Hence, qualitative robustness as defined in {(\ref{def.qualitativerobustness})} is a kind of equicontinuity concerning the weak convergence of the 
image measures of $S_{n}$ with respect to $n$.

The finite sample distribution of \RPL estimators is in general unknown. One method
to obtain approximations of this finite sample distribution is the empirical bootstrap proposed by
\cite{Efron1979,Efron1982}. As the next theorem will also contain a qualitative robustness of empirical bootstrap approximations, we need some more notation.
Recall that $\DD_n:=\frac{1}{n}\sum_{i=1}^n \delta_{(X_i,Y_i)}$ if all pairs $(X_i,Y_i)$
are independent and identically distributed with $(X_i,Y_i)\sim \P$ (abbreviation: 
$(X_i,Y_i)\stackrel{i.i.d.}{\sim}\P$).
Furthermore, we denote the distribution of the $H$-valued \RPL estimator $\fDDnLs$, by
\be 
  \Law{n}{S;\P}, \qquad n\in\N,
\ee
where $S:\PM(\cXY)\to H$ with $S(\P)=\fPLs$.
Because $\P$ is unknown but fixed, this is an unknown,  fixed probability measure
of an $H$-valued random function.
In the same manner we denote the distribution of the $H$-valued \RPL estimator $\fDDnLs$, when all pairs 
$(X_i^{(b)},Y_i^{(b)})\stackrel{i.i.d.}{\sim} \DD_n:=\frac{1}{n}\sum_{i=1}^n \delta_{(X_i,Y_i)}$, where $(X_i,Y_i)\stackrel{i.i.d.}{\sim} \P$,  by
\be 
  \Law{n}{S;\DD_n}, \quad n\in\N.
\ee
We mention that $ \Law{n}{S;\DD_n}$ denotes a distribution which can be considered itself
as a random function in an abstract sense because it depends on $\DD_n$.

We can now state our result on the qualitative robustness of 
regularized pairwise learning methods.

\begin{theorem}[\textbf{Qualitivative robustness}]\label{thm.qualitativerobust}
Let the Assumptions \ref{assumption-spaces1}, \ref{assumption-kernel1}, \ref{assumption-loss1}, and \ref{assumption-loss2} be valid.
Then, for all $\lb\in(0,\infty)$, we have:
\begin{enumerate}
\item The sequence  of \RPL estimators $(S_n)_{n\in\N}$, where $S_n:=\fDDnLs$, is qualitatively robust 
         for all Borel probability measures $\P\in\PMXY$.
\item If the metric space $\cXY$ is additionally compact, then the sequence $\Law{n}{S;\DD_n}$, $n\in\N$, of empirical bootstrap approximations
         of $\Law{n}{S;\P}$ is qualitatively robust for all Borel probability measures $\P\in\PM(\cXY)$.
\end{enumerate}
\end{theorem}

The proof of Theorem \ref{thm.qualitativerobust} is based on the 
following two results which are interesting in their own.

\begin{theorem}[\textbf{Continuity of the operator}]\label{thm.continuityofoperator}
Let the Assumptions \ref{assumption-spaces1}, \ref{assumption-kernel1}, \ref{assumption-loss1}, and \ref{assumption-loss2} be valid.
Then, for all Borel probability measures $\P\in\PMXY$
and for all $\lb\in(0,\infty)$, we have:
\begin{enumerate} 
\item The operator 
$S:  \PM(\cXY) \to H$, where $S(\P) = \fPLs$, is continuous with respect to the weak topology on $\PM(\cXY)$ and the norm topology on $H$. 
\item The operator 
$\tilde{S}:  \PM(\cXY) \to \mathcal{C}_b(\cX)$, where $S(\P) = \fPLs$, 
is continuous with respect to the weak topology on $\PM(\cXY)$ and the norm topology on $\mathcal{C}_b(\cX)$. 
\end{enumerate}
\end{theorem}

\begin{corollary}[\textbf{Continuity of the estimator}]\label{cor.continuityofestimator}
Let the Assumptions \ref{assumption-spaces1}, \ref{assumption-kernel1}, \ref{assumption-loss1}, and \ref{assumption-loss2} be valid.
For any data set $D_n\in (\cXY)^n$ denote the corresponding empirical measure by 
$\D_n:=\frac{1}{n}\sum_{i=1}^n \delta_{(x_i,y_i)}$. 
Then, for every $\lb \in (0,\infty)$ and every $n\in\N$, the mapping 
\be
S_n : \Bigl( (\cXY)^n, d_{(\cXY)^n} \Bigr) \to \bigl( H, d_H \bigr) , \quad S_n(D_n)=f_{\Ls,\D_n,\lb} ,
\ee
is continuous.
\end{corollary}

It is known that support vector machines (SVMs) are qualitatively robust for fixed values of $\lb\in(0,\infty)$ but that they \emph{can not} be qualitatively robust for the usual null-sequences $\lb_n$ needed to
obtain universal consistency, because universal consistency and
qualitatively robustness are under some mild conditions 
concurrent goals, see \cite{HableChristmann2013a} for a discussion.
Is is known that SVMs can have a somewhat weaker property called
finite-sample qualitivative robustness, see \cite{HableChristmann2011}. It is an open problem, 
whether a similar result is true for \RPL methods, and we will not address this question here.

\section{Discussion}\label{discussion}

In this paper we proved some desirable statistical robustness properties for a broad class of regularized pairwise learning methods based on kernels. Such kernel methods are used in the fields of information theoretic learning, ranking, gradient learning, and metric and similarity learning. In particular, our work complements to some respect earlier work on consistency and learning rates for minimum error entropy principles by \cite{HuFanWunZhou2013}, \cite{FanHuWuZhou2013},
\cite{HuFanWuZhou2015}, for ranking algorithms by \cite{AgarwalNiyogi2009}, for metric and similarity learning problems by \cite{CaoGuoYing2015}, and for gradient learning methods by \cite{MukherjeeZhou2006}.

The following aspects are beyond the scope of this paper and
remain open for further work.
\emph{(i)} We did not address the question of an influence function of regularized pairwise learning methods, if a bounded but \emph{non-convex} pairwise loss function is used. The main problem seems to be that in this case the function $\fPLs$ is in general not unique.
\emph{(ii)} We did not add numerical comparisons because
it is known from \cite{Principe2010} that for minimum error entropy principles such methods
can be computed in an efficient gradient descent manner. 
\emph{(iii)} It seems obvious that the results developed here for pairwise learning can in principle be established also to higher order, if one uses $U$- or $V$-statistics of 
degree $\ell>2$. E.g., if $\ell=3$,  one can consider loss functions with $9$ instead of $6$ 
arguments which yields instead of {(\ref{RPLempirical})} the optimization problem 
\be \label{RPLempiricaldegreethree}
\inf_{f\in H} \frac{1}{n^3} \sum_{i=1}^n \sum_{j=1}^n \sum_{m=1}^n L(x_i,y_i,x_j,y_j,x_m, y_m, f(x_i),f(x_j),f(x_m)) + 
\lb \hhnorm{f}\,.
\ee
We conjecture that the numerical effort to solve such problems will strongly increase with 
$\ell$.

\section{Appendix}
\subsection{Appendix A: Some Tools}
To improve the readability of the paper, we list some known results which are used in our proofs given in this  subsection.
The following theorem provides a criterion for the existence of a global minimizer.
\citet[Prop.\,6, p.\,75]{EkelandTurnbull1983} shows the existence, and 
the uniqueness is a consequence of the strict convexity.

\begin{theorem}[Existence of minimizers]\label{app1b:existence-min}
Let $E$ be a reflexive Banach space and $f:E\to \R\cup\{\infty\}$ be a convex 
and lower semi-continuous map. If there exists an $M>0$ such that
$\{x\in E: f(x)\leq M\}$ is non-empty and bounded, then 
$f$ has a global minimum, i.e., there exists an $x_{0}\in E$ with
$f(x_{0})\leq f(x)$ for all $x\in E$. 
Moreover, if $f$ is strictly convex, then $x_{0}$ is the only element minimizing $f$.
\end{theorem}

\begin{definition}[Derivatives, see e.g. {\citet[p.\,518f]{DenkowskiEtAl2003}}] 
Let $E$ and $F$ be normed spaces, 
$U\subset E$ and $V\subset F$ be open   sets,
and
$G:U\to V$ be a map. 
We say that $G$ is G\^ateaux differentiable at $x_0\in U$
if there exists a bounded linear operator $A:E\to F$ such that 
$$
\lim_{\substack{t\to 0\\t\neq 0}} \frac{\snorm{G(x_0+tx)- G(x_0) - tAx}_F}{t} = 0\, , \qquad \qquad x\in E.
$$
In this case, $A$ is called the derivative of $G$ at $x_0$, and since $A$
is uniquely determined, we write
$G'(x_0) := \frac{\partial G}{\partial E}(x_0) := A$.
Moreover, we say $G$ that Fr\'echet differentiable at $x_0$ if $A$ actually satisfies 
$$
\lim_{\substack{x\to 0\\x\neq 0}} \frac{\snorm{G(x_0+x)- G(x_0) - Ax}_F}{\snorm x_E} = 0\, .
$$
Furthermore, we say that $G$ is (G\^ateaux, Fr\'echet) differentiable if it is (G\^ateaux, Fr\'echet) differentiable at every $x_0\in U$. 
Finally, $G$ is said to be continuously differentiable if it is Fr\'echet differentiable 
and the derivative $G':U\to \ca L(E,F)$ is continuous.
\end{definition}

\begin{theorem}[Partial Fr{\'e}chet differentiability, see e.g. {\citet[Theorem 2.6, p.\,37]{Akerkar1999}}]\label{appendixAC:partialFrechet}
Let $E_1$, $E_2$, and $F$ be Banach spaces, $U_1\subset E_1$ and $U_2\subset E_2$ be   open subsets,
and $G: U_1 \times U_2 \to F$ be a continuous map.
Then $G$ is continuously differentiable  if and only if
$G$ is partially Fr{\'e}chet differentiable and the partial derivatives $\frac{\partial G}{\partial E_1}$ and $\frac{\partial G}{\partial E_2}$ 
are continuous. In this case, the 
derivative of $G$ at  $(x_1,x_2) \in U_1 \times U_2$ is given by
$$
  G'(x_1,x_2)(y_1,y_2) =\frac{\partial G}{\partial E_1}(x_1,x_2)y_1 +\frac{\partial G}{\partial E_2}(x_1,x_2)y_2\, , 
   \qquad (y_1,y_2)\in E_1\times E_2\, . 
$$
\end{theorem}

The proof of our Theorem \ref{thm-boundedderivative} heavily relies on the implicit function theorem 
in Banach spaces. 
Recall the following simplified version of this theorem, see \citet[Thm.\,4.1, Cor.\,4.2]{Akerkar1999}
Here and throughout this appendix $B_E$ denotes the open unit ball of a Banach space $E$.

\begin{theorem}[Implicit function theorem]\label{implicit-thm}
Let $E, F$ be Banach spaces and $G:E\times F\to F$ be a continuously differentiable map.
Suppose that we have $(x_0,y_0)\in E\times F$ such that $G(x_0,y_0)=0$ and
$\frac{\partial G}{\partial F}(x_0,y_0)$ is invertible. Then there exists a $\d>0$ and
a continuously differentiable map $f: x_0+\d B_E \to y_0 + \d B_F$ such that for all $x\in x_0+\d B_E$,
$y\in y_0 + \d B_F$ we have:
$G(x,y)=0$ if and only if $y=f(x)$.
Moreover, the derivative of $f$ is given by
\be \label{implicit-thm-f1}
f'(x) = -\left(\frac{\partial G}{\partial F}\bigl(x,f(x)\bigr) \right)^{-1}
  \frac{\partial G}{\partial E}\bigl(x,f(x)\bigr)\, .
\ee
\end{theorem}

%

\begin{definition}[Bochner integral]
Let $E$ be a Banach space and  $(\Om,\ca A,\mu)$ be a $\s$-finite measure space.
An $E$-valued measurable function $f:\Om\to E$ is called 
Bochner $\mu$-integrable if there exists a 
sequence $(f_n)_{n\in\N}$ of $E$-valued measurable step functions $f_n:\Om\to E$ such that 
$\lim_{n\to \infty} \int_\Om \snorm{f_n-f}_E \, d\mu = 0$.
In this case, the limit 
$\int_\Om f\,d\mu := \lim_{n\to \infty} \int_\Om f_n\, d\mu$
exists and is called the Bochner integral of $f$.
Finally, if $\mu$ is a probability measure, we sometimes write $\Ex_\mu[f]$ for this integral.
\end{definition}

\begin{theorem}[Dominated convergence theorem, see e.g. {\citet[Thm.\,3.10.12, p.\,367]{DenkowskiEtAl2003}}]\label{appendixAC:DominatedConvergenceForBochner}
Let $E$ be a Banach space, $(\O, \mathcal{A},\mu)$ be a finite measure space, 
and $(f_n)_{n\in\N}$ be a sequence of Bochner  $\mu$-integrable functions $f_n:\O\to E$.
If $\lim_{n\to\infty} \mu\{\snorm{f_n-f} \ge \e\}=0$ for every $\e>0$
and if there exists a $\mu$-integrable function $g:\Om\to \R$ with $\snorm{f_n(\omega)} \le g(\omega)$
$\mu$-almost everywhere for all $n\in\N$,
then $f$ is Bochner  $\mu$-integrable and 
$\lim_{n\to\infty} \int_\O f_n \,d\mu = \int_\O f \, d\mu$. 
\end{theorem}

%

\subsection{Appendix B: Proofs}\label{proofs}

To shorten the notation, we occationally use the abbreviations $z:=(x,y)$, $\tiz:=(\tix,\tiy)$, $Z:=(X,Y)$, $\tiZ:=(\tiX,\tiY)$,
$\cZ:=\cX\times\cY$, and $D_i L(z,\tiz,t,\tit):=D_i L(x,y,\tix,\tiy,t,\tit)$,
$D_i D_j L(z,\tiz,t,\tit):=D_i D_j L(x,y,\tix,\tiy,t,\tit)$ for $i\in\{5,6\}$ etc.

The proofs for the results given in Section \ref{loss} and Section \ref{existunique} are similar to 
corresponding results for ``classical'' loss functions of the form $L:\cX\times\cY\times\R\to[0,\infty)$ used by
support vector machines and related kernel based methods, see e.g. \citet{SC2008}.

\begin{proofof}{\myem{Proof of Lemma \ref{loss:meas-risk}}}
Since $d$ dominates the pointwise convergence, we see that, 
for fixed $x\in \cX$, the $\R$-valued map $f\mapsto f(x)$ defined on $\ca F$ is continuous with respect to $d$.
Furthermore, $\ca F\subset \sL 0 \cX$ implies that, 
for fixed $f\in \ca F$, the $\R$-valued map $x\mapsto f(x)$ defined on $\cX$ is measurable.
By a well-known result from Carath\'{e}odory, see e.g.
\citet[p.\,70]{CaVa77}, we then obtain the first assertion.
Since this implies that the maps 
$(x,y,\tix,\tiy,f) \mapsto (x,y,\tix,\tiy,f(x),f(\tix))$ 
and $(x,y,\tix,\tiy,f,f) \mapsto (x,y,\tix,\tiy,f(x),f(\tix))$ 
are measurable, we obtain the second assertion.
The third  assertion now follows from  the measurability statement 
in Tonelli-Fubini's theorem, see e.g. \citet[p.\,137]{Dudley2002}.
\qedr
\end{proofof}

\begin{proofof}{\myem{Proof of Lemma \ref{loss:convex-loss-lemma}}}
Let $c\in[0,1]$, $f,g\in \sL{0}{\cX}$, and assume that $L$ is a
convex  pairwise loss.
We immediately obtain, for all $(x,y,\tix,\tiy)\in(\cXY)^2$, 
\begin{eqnarray*}
& & L\bigl(x,y,\tix,\tiy,cf(x)+(1-c)g(x),cf(\tix)+(1-c)g(\tix)\bigr)\\
& \le  & 
       c L\bigl(x,y,\tix,\tiy,f(x),f(\tix)\bigr)
  +  (1-c) L\bigl(x,y,\tix,\tiy,g(x),g(\tix)\bigr).
\end{eqnarray*}
The linearity of integrals yields the assertion
$\RP{L}{cf+(1-c)g} \le c \RP{L}{f}+(1-c)\RP{L}{g}$. 
The case of strict convexity can be shown in an analogous manner.
\qedr
\end{proofof}

%

\begin{proofof}{\myem{Proof of Lemma \ref{loss:lipschitz-loss-lemma}}}
Because $L$ is a locally separately Lipschitz continuous pairwise loss, we have
\beqna
& & | \RPo{L}(f) - \RPo{L}(g)|  \\
& \le & \int_{(\cX\times \cY)^2} \bigl| L(x,y,\tix,\tiy,f(x),f(\tix)) -  L(x,y,\tix,\tiy,g(x),g(\tix))\bigr| 
\,d\P^2(x,y,\tix,\tiy) \\
& \le & |L|_{B,1} \int_{\cX^2} | f(x)-g(x) | +  | f(\tix) - g(\tix)| \,d\P_\cX^2(x,\tix)  \\
& \le & 2 \, |L|_{B,1} \cdot \snorm{f - g}_{L_{1}(\P_\cX)} \,,
\eeqna
which gives the assertion.
\qedr
\end{proofof}

\begin{proofof}{\myem{Proof of Lemma \ref{loss:diff-risk}}}
Because we consider $L$ as a function of its last two arguments, while the first
four arguments are held fixed, we define $L_{z,\tiz}(t,\tit):=L(z,\tiz,t,\tit)$,
where $z:=(x,y)$ and $\tiz=(\tix,\tiy)$.
We first observe that all derivatives
$(D L_{z,\tiz})(t,\tit)$ are measurable since
we assumed continuous partial derivatives.
Now let $f\in \Lx\infty{\P_{\cX}}$ and
$(f_n)_{n\in\N}\subset \Lx \infty{\P_{\cX}}$ be a sequence 
with $f_{n}\neq 0$, $n\geq 1$, 
and $\lim_{n\to \infty}\inorm{f_{n}}= 0$.
Without loss of generality, we additionally assume for later use 
that $\inorm{f_{n}}\leq 1$ for all $n\geq 1$. 
For $z=(x,y), \tiz=(\tix,\tiy)\in \cZ$ and $n\geq 1$, we now define
\begin{eqnarray*}
G_n(z,\tiz) & := &  
\sqrt{2} \cdot \Bigl| \Lzz\bigl(f(x)+f_n(x), f(\tix)+f_n(\tix)\bigr) 
       - \Lzz\bigl(f(x),f(\tix)\bigr) \\
& &      \,\,\,\,\,\,\,\, - \bigl\langle (D\Lzz)(f(x),f(\tix)), 
             \bigl(f_n(x),f_n(\tix)\bigr) \bigr\rangle \Bigr| 
             \Bigl/ \snorm{(f_n(x),f_n(\tix))}_2 \,,
\end{eqnarray*}
if $\snorm{(f_n(x),f_n(\tix))}_2 \ne 0$, and $G_n(z,\tiz):=0$ otherwise.
We obtain
\begin{eqnarray}\nonumber
&&
\biggl| \frac{\RP L{f+f_{n}} - \RP L f - \RPs{L}{f}{f_n}}
             {\inorm{f_{n}}}   \biggl| \\ \nonumber 
& \le &
\int_{(\cXY)^2}  \frac{1}{\inorm{f_n}} \, \cdot\,    
\Bigl| \Lzz\bigl(f(x)+f_n(x),f(\tix)+f_n(\tix)\bigr) 
         - \Lzz\bigl(f(x),f(\tix)\bigr) \\ \nonumber
& & ~~~~~~~~~~~~~~~~~~~  - \bigl\langle (D\Lzz)(f(x),f(\tix)), 
        \bigl(f_n(x),f_n(\tix)\bigr) \bigr\rangle \Bigr| \, d\P^2(x,y,\tix,\tiy)\nonumber \\
& \le &
\int_{(\cXY)^2}  \frac{\sqrt{2}}{\snorm{(f_n(x),f_n(\tix))}_2} \, \cdot\,  \label{loss:diff-risk-h1}  \\
& & \,\,\,\,\,\, \Bigl| \Lzz\bigl(f(x)+f_n(x),f(\tix)+f_n(\tix)\bigr)  
                         - \Lzz\bigl(f(x),f(\tix)\bigr)   \nonumber \\       
& &\,\,\,\,\,\,\,\,\,  - \bigl\langle (D\Lzz)(f(x),f(\tix)),
 \bigl(f_n(x),f_n(\tix)\bigr) \bigr\rangle \Bigr| \, d\P^2(x,y,\tix,\tiy) 
\nonumber  \\
& = &
\int_{\cZ^2}    
G_n(z,\tiz)\,  d\P^2(z,\tiz) \label{loss:diff-risk-h1c}
\end{eqnarray}
for all $n\geq 1$, where the well-known relationship 
$\snorm{v}_2 \le \sqrt{2}\,\inorm{v}$, where $v\in\R^2$,
was used in {(\ref{loss:diff-risk-h1})}.
Furthermore, for $z,\tiz\in\cZ$, the definition of $G_n$ and
the definition of the (total) derivative
$D\Lzz$ obviously yield
\begin{equation}\label{loss:diff-risk-h2}
\lim_{n\to \infty} G_n(z,\tiz) = 0.
\end{equation}
Denote the gradient of $\Lzz$ by $\nabla\Lzz$.
For $z,\tiz\in \cZ$ and $n\ge 1$ with $f_{n}(x)\ne 0$,
the mean value theorem for functions from $\R^2$ to $\R$ shows 
that there exists some $g_n(z,\tiz)$ with
\begin{eqnarray}\label{loss:diff-risk-h2b}
& & \Lzz\bigl(f(x)+f_n(x), f(\tix)+f_n(\tix)\bigr) 
       - \Lzz\bigl(f(x),f(\tix)\bigr) \\
& = &  \Bigl\langle (\nabla\Lzz)\bigl((1-g_n(z,\tiz))f(x)+g_n(z,\tiz)f_n(x), \nonumber \\
& &   ~~~~~~~~~~~~~ (1-g_n(z,\tiz))f(\tix)+g_n(z,\tiz)f_n(\tix)\bigr), ~~
             \bigl(f_n(x),f_n(\tix)\bigr) \Bigr\rangle \,. \nonumber
\end{eqnarray}
Because $L$ has by assumption uniformly bounded partial 
derivatives $D_5 L$ and $D_6 L$, to be more precise, there exists a constant 
$c_L\in[0,\infty)$ such that
$$
\sup_{x,\tix\in\cX,~y,\tiy\in\cY,~t,\tit\in\R^2} 
\bigl|D_i L(x,y,\tix,\tiy,t,\tit)\bigr|
\le c_L\,, \qquad i\in\{5,6\}.
$$
If we combine this with the equality in {(\ref{loss:diff-risk-h2b})},
we obtain
\beq \label{loss:diff-risk-h4}
& & \Bigl| \Lzz\bigl(f(x)+f_n(x), f(\tix)+f_n(\tix)\bigr) 
       - \Lzz\bigl(f(x),f(\tix)\bigr) \Bigr| \\
& \le & 2 c_L \bigl( |f_n(x)| + |f_n(\tix)| \bigr)
        ~ \le ~ 4 c_L \inorm{f_n} ~ \longrightarrow 0, \quad n\to\infty\, . \nonumber
\eeq
Combining {(\ref{loss:diff-risk-h1c})},
{(\ref{loss:diff-risk-h2b})}, and {(\ref{loss:diff-risk-h4})}, 
we get 
$\int_{\cZ^2} G_n(z,\tiz) \, d\,\P^2(z,\tiz) < \infty$.
Hence, $G_n$ is a non-negative convergent dominating function and the
assertion follows from Lebesgue's theorem.
\qedr
\end{proofof}


\begin{proofof}{\textbf{Proof of Lemma \ref{weakcong}}} 
Since the closed ball $B(f_0, r)$ of the Hilbert space $H$ is weakly compact, there exists a 
subsequence $(f_{j_\ell})_{\ell=1}^\infty$ weakly converging to some $f^{} \in B(f_0, r)$. That is, 
\begin{equation}\label{weakc} 
\lim_{\ell\to\infty} \langle
f_{j_\ell}, f\rangle_H = \langle f^{*}, f\rangle_H,
\qquad \forall f\in H. 
\end{equation}

Let $f= f^{*}$ in (\ref{weakc}). Then 
$$ \|f^{*}\|_H^2 = \langle f^{*}, f^{*}\rangle_H = \lim_{\ell\to\infty} \langle
f_{j_\ell}, f^{*}\rangle_H \leq \underline{\lim}_{\ell\to\infty} \|f_{j_\ell}\|_H \|f^{*}\|_H $$
which implies $\|f^{*}\|_H \leq \underline{\lim}_{\ell\to\infty} \|f_{j_\ell}\|_H$. 

Let $f=\Phi(x)$, $x\in\cX$, in (\ref{weakc}). Then the reproducing property of the kernel yields 
$$ \lim_{\ell\to \infty} f_{j_\ell} (x)= \lim_{\ell\to\infty} \langle
f_{j_\ell}, \Phi(x)\rangle_H = \langle f^{*}, \Phi(x)\rangle_H = f^{*} (x). $$
This is true for any $x\in\cX.$ So (\ref{converge}) is verified. 
\qedr
\end{proofof} 

\begin{proofof}{\textbf{Proof of Theorem \ref{thmexist}}}
For every $\ell \in\NN$, we take a function $f_\ell\in H$ such that 
\begin{equation}\label{deffk}
 {\mathcal R}_{L, \P}(f_\ell) + \lambda \|f_\ell\|_H^2 \leq \inf_{f\in H} {\mathcal R}_{L, \P}(f) + \lambda \|f\|_H^2 + \frac{1}{\ell}. 
 \end{equation}
Taking $f=f_0$, we find that 
$$ \lambda \|f_\ell\|_H^2 \leq {\mathcal R}_{L, \P}(f_0) + \lb \hhnorm{f_0} + 1 $$
and 
$$ f_\ell \in B(0, r), \ \hbox{with} \ r= \frac{{\mathcal R}_{L, \P}(f_0) + \lb \hhnorm{f_0}+ 1}{\lambda}. $$
Now we apply Lemma \ref{weakcong}. We know that there exists a subsequence $(f_{j_\ell})_{\ell=1}^\infty$ and $f^{*} \in B(0, r)$ such that
 $\|f^{*}\|_H \leq \underline{\lim}_{\ell\to\infty} \|f_{j_\ell}\|_H$ and (\ref{converge}) is valid. 
 
Consider $ {\mathcal R}_{L, \P}(f_\ell)$. By the Lipschitz continuity of $L$, the integrated function is bounded by
$$ 
L\left(x, y, \tilde{x}, \tilde{y}, f_{\ell}(x), f_{\ell}(\tilde{x})\right) 
\leq 
L\left(x, y, \tilde{x}, \tilde{y}, f_0(x), f_0(\tix)\right) + 4 |L|_1 \inorm{k} r. 
$$
The upper bound is integrable with respect to $\P^2$. Also, for any $(x, y, \tilde{x}, \tilde{y})\in(\cXY)^2$, by the continuity of $L$ and (\ref{converge}), we have
$$ 
\lim_{\ell\to \infty} L\left(x, y, \tilde{x}, \tilde{y}, f_{\ell}(x), f_{\ell}(\tilde{x})\right) = L\left(x, y, \tilde{x}, \tilde{y}, f^{*}(x), f^{*}(\tilde{x})\right). 
$$
So by the Lebesgue Dominated Theorem, we have 
$$
\lim_{\ell\to \infty} {\mathcal R}_{L, \P}(f_\ell) = {\mathcal R}_{L, \P}(f^{*}). $$ 
Then we take $\underline{\lim}$ on both sides of (\ref{deffk}) and find from $\|f^{*}\|_H \leq \underline{\lim}_{\ell\to\infty} \|f_{j_\ell}\|_H$ that 
$$ {\mathcal R}_{L, \P}(f^{*}) + \lambda \|f^{*}\|_H \leq \inf_{f\in H} {\mathcal R}_{L, \P}(f) + \lambda \|f\|_H^2. $$ 
It means that $f^{*}$ is a minimizer $f_{L, \P, \lambda}$. This proves our statement. 
\qedr
\end{proofof}

\begin{proofof}{\myem{Proof of Lemma \ref{infinite:fp-unique}}}
Let us assume that the map $f\mapsto \RP{L}{f} + \lb \hhnorm{f} $ has two minimizers $f_1,f_2\in H$
with $f_1\neq f_2$.
Recall that the parallelogram law 
$\snorm{x+x'}^2+\snorm{x-x'}^2=2\snorm{x}^2+2\snorm{x'}^2$
is valid for all points $x$ and $x'$ in a Hilbert space, \emph{cf.} \citet[Thm.\,3.7.7, p.\,310]{DenkowskiEtAl2003}.
Therefore, we have 
$\hhnorm{\frac{1}{2}(f_1+f_2)} <  \frac{1}{2} \hhnorm{f_1} + \frac{1}{2}\hhnorm{f_2}$.
As $L$ is a convex pairwise loss, the map 
$f\mapsto \RP{L}{f}$ is convex by Lemma \ref{loss:convex-loss-lemma}.
This together with 
$\RP L {f_1} + \lb \hhnorm{f_1} = \RP{L}{f_2} + \lb \hhnorm{f_2}$
then shows for $f^*:= \frac{1}{2}(f_1+f_2)$ that
$$
\RP{L}{f^*} + \lb\hhnorm{f^*}  < \RP{L}{f_1} + \lb \hhnorm{f_1} \, ,
$$
i.e., $f_1$ is not a minimizer of $f\mapsto \RP{L}{f} + \lb \hhnorm{f}$.
Consequently, the assumption that there are two  minimizers is false.
\qedr
\end{proofof}

\begin{proofof}{\myem{Proof of Theorem \ref{infinite:fp-exist}}}
Since the kernel $k$ of $H$ is measurable, $H$ consists of measurable functions, see
e.g. \citet[Lem.\,4.24]{SC2008}.
Moreover, $k$ is bounded and thus $\id: H\to \Lx\infty {\P_{\cX}}$ is continuous, see
e.g. \citet[Lem.\,4.23]{SC2008}.
In addition, we have $L(z,\tiz,t,\tit) < \infty$ for all 
$(z,\tiz,t,\tit)\in \cZ^2\times \R^2$.
Recall that every convex function $g:\R^2\to\R\cup\{\infty\}$, which is 
not identically equal to $+\infty$, 
is continuous on the interior of its effective domain $\Dom{g}:=\{t\in\R; g(t)<\infty\}$, 
see e.g. \citet[Cor.\,2.3 on p.\,12]{EkelandTemam1999}.
Hence $L$ is a continuous pairwise loss by the convexity of $L$.
Therefore, Lemma \ref{loss:lipschitz-loss-lemma} shows that 
$\RPo L :\Lx\infty {\P_{\cX}}\to \R$ is continuous, and 
hence
$\RPo L :H\to \R$ is continuous. In addition, Lemma \ref{loss:convex-loss-lemma} 
provides the convexity of this map.
Furthermore, $f\mapsto \lb \hhnorm{f}$ is also  convex and continuous, 
which yields the continuity and the convexity of the map 
$f\mapsto \cR_{L,\P,\lb}^{reg}(f)$.
Now consider the set 
$A:=\bigl\{f\in H: \cR_{L,\P,\lb}^{reg}(f) \le \RP{L}{0} \bigr\}$.
We obviously have $0\in A$. 
In addition, $f\in A$ implies $\lb \hhnorm{f} \le \RP{L}{0}$, 
and hence
$A\subset  (\RP{L}{0}/\lb)^{1/2} \bar{B}_{H}$, 
where $\bar{B}_H$ denotes the closed unit ball of $H$. 
Hence $A$ is a non-empty and bounded subset of $H$ 
and thus Theorem \ref{app1b:existence-min} gives the 
existence of a minimizer $\fP$.
\qedr
\end{proofof}

\begin{proofof}{\myem{Proof of Corollary \ref{infinite:exist-and-unique-mar-based}}}
Because $L$ is a separately Lipschitz continuous pairwise loss, we have,
for all $(z,\tiz,t,\tit)\in \cZ^2\times\R^2$,  
\beqna
L(z,\tiz,t,\tit)  & = & 
  L(z,\tiz,0,0) + L(z,\tiz,t,\tit) - L(z,\tiz,0,0)  \\
& \le &   L(z,\tiz,0,0) + 2 |L|_1\, \bigl( |t| + |\tit| \bigr).
\eeqna
The assumption $\RP{L}{0} < \infty$ yields $\cR_{L,\P}(f) \le \cR_{L,\P}(0)+4 |L|_1 \snorm{f}_{L_1(\P_\cX)} < \infty$. 
As $L$ is a convex pairwise loss, 
Lemma \ref{infinite:fp-unique} and Theorem \ref{infinite:fp-exist}
yield the existence and the uniqueness of $\fP$
and {(\ref{naivefPbound})} equals {(\ref{infinite-fp-bound})}.
\qedr 
\end{proofof}


\begin{proofof}{\myem{Proof of Lemma \ref{L*prop}}}
Follows immediately from the definition of $\Ls$. \qed 
\end{proofof}

\begin{proofof}{\myem{Proof of Lemma \ref{props}}}
\emph{(i)} Obviously, we have 
$$
  \inf_{f\in\mathcal{L}_0(\cX)} \Ls(x,y,\tix,\tiy,f(x),f(\tix)) 
  \le 
  \Ls(x,y,\tix,\tiy,0,0)=0.
$$
\emph{(ii)} We have for all $f \in H$ that
\begin{eqnarray*}
|\RP{\Ls}{f}| & = &  |\Ex_{\P^2} L(X,Y,\tiX,\tiY,f(X),f(\tiX)) - L(X,Y,\tiX,\tiY,0,0)| \\
& \le & \Ex_{\P^2} |L(X,Y,\tiX,\tiY,f(X),f(\tiX)) - L(X,Y,\tiX,\tiY,0,0)| \\
& \le & |L|_1 \, \Ex_{\P^2} \bigl(|f(X)| +  |f(\tiX)|\bigr) \\
& \le &  2 |L|_1 \, \Ex_{\P_\cX} |f(X)|,
\end{eqnarray*}
which proves (\ref{riskbound}). Equation (\ref{riskbound2}) follows
from $\RPr{\Ls}{f} = \RP{\Ls}{f} + \lambda \hhnorm{f}$.\\
\emph{(iii)} As $0 \in H$, we obtain
$\inf_{f \in H} \RPr{\Ls}{f} \le \RPr{\Ls}{0} = 0$ and
the same reasoning holds for $\inf_{f \in H} \RP{\Ls}{f}$.\\
\emph{(iv)} Due to \emph{(iii)} we have $\RPr{\Ls}{\fPLs} \le 0$.
As $L\ge 0$ we obtain
\begin{eqnarray*}
  \lambda \hhnorm{\fPLs}
  & \le & - \RP{\Ls}{\fPLs} \\
  & = &  \Ex_{\P^2} \bigl( L(X,Y,\tiX,\tiY,0,0) - L(X,Y,\tiX,\tiY,\fPLs(X),\fPLs(\tiX)) \bigr) \\
  & \le & \Ex_{\P^2} L(X,Y,\tiX,\tiY,0,0) = \RP{L}{0}.
\end{eqnarray*}
Using similar arguments as above, we obtain
\begin{eqnarray*}
0 & \le & - \RPr{\Ls}{\fPLs} \\
  & = & \Ex_{\P^2} \bigl(L(X,Y,\tiX,\tiY,0,0) - L(X,Y,\tiX,\tiY,\fPLs(X),\fPLs(\tiX))\bigr) - \lambda \hhnorm{\fPLs} \\
& \le & \Ex_{\P^2} L(X,Y,\tiX,\tiY,0,0) = \RP{L}{0}.
\end{eqnarray*}
Furthermore, we obtain
\begin{eqnarray*}
  - 2 |L|_1 \, \Ex_{\P_\cX} |\fPLs(X)| + \lb \hhnorm{\fPLs}
  & \le & \RPr{\Ls}{\fPLs} 
   \le  \RPr{\Ls}{0} = 0.
\end{eqnarray*}
This yields (\ref{fbound1}). Using {(\ref{kernel-prop1})}, {(\ref{kernel-prop3})}, and
(\ref{fbound1}), we obtain for $\fPLs \ne 0$ that
\begin{eqnarray*}
\inorm{\fPLs} & \le & \inorm{k} \hnorm{\fPLs} \\
   & \le & \inorm{k} \sqrt{(2/ \lb) \, |L|_1 \, \Ex_{\P_\cX}|\fPLs(X)|} \\
          & \le & \inorm{k} \sqrt{(2/ \lb) \, |L|_1 \inorm{\fPLs}} < \infty \,.
\end{eqnarray*}
Hence $\inorm{\fPLs} \le \frac{2}{\lb} \, \inorm{k}^2 \, |L|_1$.
The case $\fPLs = 0$ is trivial.

\emph{(v)} This follows immediately from the definition of $\Ls$, because we just subtract a term, which is constant w.r.t.
the last two arguments of $\Ls$.
\qedr 
\end{proofof}

\begin{proofof}{\myem{Proof of Lemma \ref{proprisk}}}
The inequality 
$|\RP{\Ls}{f}| \le 2\,|L|_1 \Ex_{\P_\cX} |f(X)| < \infty$
for $f \in L_1(\P_\cX)$ follows from (\ref{riskbound}).
Then (\ref{riskbound2}) yields
$\RPr{\Ls}{f} \ge  -|L|_1 \Ex_{\P_\cX} |f(X)| + \lambda \hnorm{f}^2
> -\infty$.
\qedr
\end{proofof}

\begin{proofof}{\myem{Proof of Lemma \ref{convrisk2}}}
Lemma \ref{L*prop} yields that $\Ls$ is (strictly)
convex. Trivially $\RPo{\Ls}$ is also (strictly) convex. Further 
$f \mapsto \lambda \hhnorm{f}$ is strictly convex, and hence the map 
$f \mapsto \RPr{\Ls}{f}=\RP{\Ls}{f}+\lambda \hhnorm{f}$ is strictly convex. 
\qedr
\end{proofof}

\begin{proofof}{\myem{Proof of Theorem \ref{infinite:fp-unique2}}}
The proof of part (i) is almost identical to the proof of Lemma \ref{infinite:fp-unique}.
We only have to use Lemma \ref{convrisk2}  instead of Lemma \ref{loss:convex-loss-lemma}.
\emph{(ii)} This condition implies that
$|\RP{\Ls}{f}| < \infty$, see Lemma \ref{proprisk},
and the assertion follows from \emph{(i)}.
\qedr
\end{proofof}

\begin{proofof}{\myem{Proof of Theorem \ref{infinite:fp-exist2}}}
Because the  proof is very similar to the proof of Theorem \ref{infinite:fp-exist}, we omit it.
We also refer to the proof of Theorem 6 by \citet{ChristmannVanMessemSteinwart2009} for details. 
The uniqueness of $\fLsP$ follows immediately from Theorem \ref{infinite:fp-unique2}(ii), because 
the boundedness of $k$ guarantees $\inorm{f}\le \inorm{k}\,\hnorm{f}$, see
{(\ref{kernel-prop1})}. 
\qedr
\end{proofof}


Before we can prove Theorem \ref{thm-boundedderivative},
we need the following results.

\begin{lemma} \label{Lemma-root}
Let the Assumptions \ref{assumption-spaces1}, \ref{assumption-kernel1}, 
and \ref{assumption-loss1} be satisfied. Let $\fPLs\in H$ be any fixed minimizer of $\min_{f\in H} \bigl( \RP{\Ls}{f}+\lb \hhnorm{f}\bigr)$.
Then we have, for any $g\in H$,
\begin{eqnarray*}
& & \Ex_{\P^2} \biggl[ \DfiveL{\fPLs} g (X) \\
& &  ~~+ \DsixL{\fPLs} g(\tiX) \biggr] 
+ 2 \lambda  \langle \fPLs, g\rangle_H = 0.
\end{eqnarray*}
\end{lemma}

\begin{proofof}{\textbf{Proof of Lemma \ref{Lemma-root}}}
Because $\Ls$ and $\lb$ are fixed, we write $f_\P:=\fPLs$
to shorten the notation in the proof. 
Let $g\in H$. We define a continuous function 
\be \label{defGtilde}
 \tilde{G}: [-1, 1] \to \RR, \quad  \tilde{G}(t) = {\mathcal R}_{\Ls, \P}(f_{\P} + t g) + \lambda \hhnorm{f_{\P} + t g}. 
\ee
Recall that the partial derivatives of $L$ and $\Ls$ w.r.t. to the last two arguments are identical because 
$L$ and $\Ls$ differ only by the term $L(x,y,\tix,\tiy,0,0)$.
Observe that for $t\not= 0$,
\begin{eqnarray*}
\frac{\tilde{G}(t) - \tilde{G}(0)}{t} &=& \int_{(\cXY)^2} \frac{1}{t}\biggl(L\left(x, y, \tilde{x}, \tilde{y}, f_\P(x) + t g(x), f_\P(\tilde{x}) + t g(\tilde{x})\right) \\
&& \qquad \quad \qquad - L\left(x, y, \tilde{x}, \tilde{y}, f_\P(x), f_\P(\tilde{x})\right)\biggr) d\P^2(x, y, \tilde{x}, \tilde{y}) \\
&&  + 2 \lambda  \langle f_\P, g\rangle_H + t \|g\|_H^2.
\end{eqnarray*}
By the separate Lipschitz continuity of $L$, the absolute value of the above integrand is bounded by
$$ |L|_1 \left(|g(x)| + |g(\tilde{x})|\right) \leq 2 |L|_1 \|g\|_\infty < \infty\,. $$
Also, for any $(x, y), (\tilde{x}, \tilde{y})\in\cXY$, we have
\begin{eqnarray*}
 \lim_{t\to 0} \frac{1}{t}\biggl( L\big(x, y, \tilde{x}, \tilde{y}, f_\P(x) + t g(x), f_\P(\tilde{x}) + t g(\tilde{x})\big)  
    - L\big(x, y, \tilde{x}, \tilde{y}, f_\P(x), f_\P(\tilde{x})\big)\biggr) \\
=   \DfiveLsmall{f_{\P}} g(x) +   \DsixLsmall{f_{\P}} g(\tilde{x}).
\end{eqnarray*}
An application of  Lebesgue's dominated convergence theorem yields 
\begin{eqnarray*}
\lim_{t\to 0} \frac{\tilde{G}(t) - \tilde{G}(0)}{t} &=& \int_{(\cXY)^2} \biggl(\DfiveLsmall{f_{\P}} g(x) \\
&& \qquad + \DsixLsmall{f_{\P}} g(\tilde{x}) \biggr) d\P^2(\xyxy) \\
&& + 2 \lambda  \langle f_\P, g\rangle_H.
\end{eqnarray*}
Since $\tilde{G}(t) - \tilde{G}(0) \geq 0$ for any $t$ by definition of $\tilde{G}$, we know that the term on the right hand of the above equality is greater or equal to $0$.
This inequality is also true for the function $-g$. So the desired identity follows.
\qedr 
\end{proofof}

\begin{definition}\label{modcon}
We define the local modulus of continuity for the second order derivatives of a pairwise loss function $L$ with respect to the last two variables as 
\begin{eqnarray}
&&\omega (h)_r := \sup \biggl\{\left|D_i D_j L\bigl(x, y, \tilde{x}, \tilde{y}, f, \tilde{f}\bigr) - D_i D_j L\bigl(x, y, \tilde{x}, \tilde{y}, g, \tilde{g}\bigr) \right|: \ x, \tilde{x} \in {\mathcal X}, \nonumber \\
&& \quad y, \tilde{y}\in {\mathcal Y},  f, \tilde{f}, g, \tilde{g} \in [-r, r], |f-g| \leq h, |\tilde{f} - \tilde{g}| \leq h, i, j \in \{5, 6\}\biggr\}  \label{OmegaDef}
\end{eqnarray}
where $h, r>0$. 
\end{definition}

If the sets ${\mathcal X}$ and ${\mathcal Y}$ are bounded, the continuity of the second order derivatives of $L$ 
implies that $\lim_{h\to 0} \omega (h)_r =0$ uniformly with respect to $r >0$. 

Let $\P,\Q\in \mathcal{M}((\cXY)^2)$ and $\varepsilon\in\R$.
Define the signed measure $\P_\varepsilon:=(1-\varepsilon)\P+\varepsilon \Q$. 
Note that $\P_\varepsilon$ is a probability distribution if $\varepsilon\in[0,1]$.

The key property we need to prove Theorem \ref{thm-boundedderivative} is formulated in the following result.

\begin{theorem}\label{THM.PropertiesOfG}
Let $L$ satisfy the Assumptions \ref{assumption-spaces1}, \ref{assumption-kernel1}, \ref{assumption-loss1}, and \ref{assumption-loss2}.  
If $\lim_{h\to 0} \omega (h)_r =0$ for any fixed $r>0$, then the function $G: \RR \times H \to H$ defined by 
\begin{eqnarray*}
G(\varepsilon, f) & = & 2 \lambda f + \E_{\P_\varepsilon^2} \Bigl[D_5L(X,Y, \tiX,\tiY, f(X), f(\tilde{X})) \Phi(X) \\
& & \qquad \qquad \quad + D_6 L(X,Y,\tiX,\tiY, f(X), f(\tilde{X})) \Phi(\tilde{X})\Bigr]
\end{eqnarray*}
is continuously differentiable. Moreover, $\frac{\partial G}{\partial H}(0, f)$ is invertible for any $f\in H$. 
\end{theorem}

\begin{proofof}{\textbf{Proof of Theorem \ref{THM.PropertiesOfG}}}
By Theorem \ref{appendixAC:partialFrechet}, we only need to show that the partial derivatives 
$\frac{\partial G}{\partial \varepsilon}$ and $\frac{\partial G}{\partial H}$ are continuous. 

To shorten the notation in the proof, we denote the random functions
$\DfiveL{f}$ by $D_5 L_f$ and $\DsixL{f}$ by $D_6 L_f$, respectively. Analogously we denote second order
partial derivatives of $L$ by $D_i D_j L_f$, where $i,j\in\{5,6\}$.

Note that for $\varepsilon \in \RR$ and $f\in H$, 
\begin{eqnarray}
\frac{\partial G}{\partial \varepsilon}(\varepsilon, f)
 & = &   -2 (1-\varepsilon) \Ex_{\P\otimes \P} \bigl( D_5 L_f \Phi(X) + D_6 L_f \Phi(\tiX) \bigr)  \label{dGdeps}\\
&     & + (1-2\varepsilon) \Ex_{\P\otimes \Q} \bigl( D_5 L_f \Phi(X) + D_6 L_f \Phi(\tiX) \bigr)  \nonumber \\
&     & + (1-2\varepsilon) \Ex_{\Q\otimes \P} \bigl( D_5 L_f \Phi(X) + D_6 L_f \Phi(\tiX) \bigr)  \nonumber \\
&     & + 2\varepsilon \Ex_{\Q\otimes \Q} \bigl( D_5 L_f \Phi(X) + D_6 L_f \Phi(\tiX) \bigr). \nonumber
\end{eqnarray}
Then for $\varepsilon, \tilde{\varepsilon} \in \RR$ and $f, \tilde{f}\in H$, we have 
\begin{eqnarray*}
\frac{\partial G}{\partial \varepsilon}(\varepsilon, f) -\frac{\partial G}{\partial \varepsilon}(\tilde{\varepsilon}, \tilde{f}) 
& = &
\left(\frac{\partial G}{\partial \varepsilon}(\varepsilon, f) -\frac{\partial G}{\partial \varepsilon}(\varepsilon, \tilde{f})\right)
 + 
 \left(\frac{\partial G}{\partial \varepsilon}(\varepsilon, \tilde{f}) -\frac{\partial G}{\partial \varepsilon}(\tilde{\varepsilon}, \tilde{f})\right) \\
&  =:  &
 \partial G_1 + \partial G_2. 
\end{eqnarray*}
Here 
\begin{eqnarray*}
\partial G_1 
& = & -2(1-\varepsilon) \Ex_{\P\otimes \P} \left[\bigl( D_5 L_f - D_5 L_{\tilde{f}} \bigr) \Phi(X) +  \bigl( D_6 L_f - D_6 L_{\tilde{f}} \bigr) \Phi(\tiX)\right] \\
&    & + (1-2\varepsilon) \Ex_{\P\otimes \Q} \left[\bigl( D_5 L_f - D_5 L_{\tilde{f}} \bigr) \Phi(X) +  \bigl( D_6 L_f - D_6 L_{\tilde{f}} \bigr) \Phi(\tiX) \right]\\
&    & + (1-2\varepsilon) \Ex_{\Q\otimes \P} \left[ \bigl( D_5 L_f - D_5 L_{\tilde{f}} \bigr) \Phi(X) +  \bigl( D_6 L_f - D_6 L_{\tilde{f}} \bigr) \Phi(\tiX) \right]\\
&    & + 2\varepsilon \Ex_{\Q\otimes \Q} \left[ \bigl( D_5 L_f - D_5 L_{\tilde{f}} \bigr) \Phi(X) +  \bigl( D_6 L_f - D_6 L_{\tilde{f}} \bigr) \Phi(\tiX) \right]. 
\end{eqnarray*}

Applying {(\ref{kernel-prop2})} and {(\ref{loss-assump2})} yield 
$$ 
\|\partial G_1\|_H 
\leq 
\left( 2 |1-\varepsilon| + 2 |1-2\varepsilon| + 2|\varepsilon|\right)  \cdot 2 \cdot\left( \inorm{k} 
\cdot 2 \cdot c_{L, 2} \inorm{k}  \|f-\tilde{f}\|_H\right).$$
Moreover, 
\begin{eqnarray*} 
\partial G_2 
& = &  ~~2(\varepsilon -\tilde{\varepsilon}) \Ex_{\P\otimes \P}  \left[ D_5 L_{\tilde{f}} \Phi(X) +  D_6 L_{\tilde{f}} \Phi(\tiX)\right] \\
&    & + 2 (\tilde{\varepsilon}-\varepsilon) \Ex_{\P\otimes \Q} \left[ D_5 L_{\tilde{f}} \Phi(X) +  D_6 L_{\tilde{f}} \Phi(\tiX)\right] \\
&    & + 2 (\tilde{\varepsilon}-\varepsilon) \Ex_{\Q\otimes \P} \left[ D_5 L_{\tilde{f}} \Phi(X) +  D_6 L_{\tilde{f}} \Phi(\tiX)\right] \\
&    & + 2(\varepsilon - \tilde{\varepsilon}) \Ex_{\Q\otimes \Q} \left[ D_5 L_{\tilde{f}} \Phi(X) +  D_6 L_{\tilde{f}} \Phi(\tiX)\right] .
\end{eqnarray*}
Hence by {(\ref{loss-assump1})} we have 
$$
 \|\partial G_2\|_H \leq 8 \cdot 2 \cdot |\varepsilon -\tilde{\varepsilon}|  c_{L, 1} \inorm{k}  .
 $$
Thus 
$$ 
\left \| \frac{\partial G}{\partial \varepsilon}(\varepsilon, f) 
-
\frac{\partial G}{\partial \varepsilon}(\tilde{\varepsilon}, \tilde{f})\right\|_H 
\leq 
\left(4 + 8 |\varepsilon| \right) 4 \inorm{k}^2 c_{L, 2} \|f-\tilde{f}\|_H + 16 c_{L, 1} \inorm{k}  |\varepsilon -\tilde{\varepsilon}|. $$
This proves the continuity of the partial derivative $\frac{\partial G}{\partial \varepsilon}$. 

The other partial derivative $\frac{\partial G}{\partial H}$ can be expressed with $\varepsilon \in \RR$ and $f\in H$ as 
\begin{eqnarray*}
\frac{\partial G}{\partial H}(\varepsilon, f) 
& = & 2 \lambda id_H + \Ex_{\P_\varepsilon \otimes \P_\varepsilon}
\biggl[D_5 D_5 L_f \Phi(X) \otimes \Phi(X) + D_6 D_5 L_f \Phi(X) \otimes \Phi(\tilde{X}) \\
&& + D_5 D_6 L_f \Phi(\tilde{X}) \otimes \Phi(X) + D_6 D_6 L_f \Phi(\tilde{X}) \otimes \Phi(\tilde{X}) \biggr].
\end{eqnarray*}

To prove its continuity, we consider first the following difference with $\tilde{f}\in H$ 
\begin{eqnarray*}
&  & \frac{\partial G}{\partial H}(\varepsilon, f) - \frac{\partial G}{\partial H}(\varepsilon, \tilde{f}) \\
& = & \Ex_{\P_\varepsilon \otimes \P_\varepsilon}
\biggl[D_5 \left(D_5 L_f - D_5 L_{\tilde{f}} \right) \Phi(X) \otimes \Phi(X) + D_6 \left( D_5 L_f - D_5 L_{\tilde{f}} \right) \Phi(X) \otimes \Phi(\tilde{X}) \\
& & \qquad + D_5 \left( D_6 L_f  - D_6 L_{\tilde{f}} \right) \Phi(\tilde{X}) \otimes \Phi(X) 
       + D_6 \left( D_6 L_f - D_6 L_{\tilde{f}} \right)\Phi(\tilde{X}) \otimes \Phi(\tilde{X}) \biggr].
\end{eqnarray*}
By the definition of the local modulus of continuity for the second order derivatives of $L$, 
see Definition \ref{modcon}, 
and the bound $\|\Phi(\tilde{X}) \otimes \Phi(\tilde{X})\|_{{\mathcal L}(H, H)}\leq \inorm{k}^2$, if  $f, \tilde{f} \in \{g\in H: \|g\|_H \leq r\}$ ,
we obtain the bound  
$$ 
\left\|\frac{\partial G}{\partial H}(\varepsilon, f) 
- 
\frac{\partial G}{\partial H}(\varepsilon, \tilde{f})\right\|_{{\mathcal L}(H, H)} 
\leq 
4 \inorm{k}^2 \omega (\inorm{k} \|f- \tilde{f}\|_H)_{r \inorm{k} }. 
$$

The second difference we need to consider is the following sum of four terms, where the integrands are the same but the factors and the probability measures differ:
\begin{eqnarray*}
& & \frac{\partial G}{\partial H}(\varepsilon, \tilde{f}) - \frac{\partial G}{\partial H}(\tilde{\varepsilon}, \tilde{f}) \\
& = & (\tilde{\varepsilon}-\varepsilon)(2- \tilde{\varepsilon}-\varepsilon) \Ex_{\P\otimes \P} 
\biggl[D_5 D_5 L_f  \Phi(X) \otimes \Phi(X) + D_6 D_5 L_f \Phi(X) \otimes \Phi(\tilde{X}) \\
& & \qquad + D_5 D_6 L_f \Phi(\tilde{X}) \otimes \Phi(X) + D_6 D_6 L_f \Phi(\tilde{X}) \otimes \Phi(\tilde{X}) \biggr] \\
& & + (\varepsilon - \tilde{\varepsilon})(1- \tilde{\varepsilon}-\varepsilon) \Ex_{\P\otimes \Q}
 \biggl[D_5 D_5 L_f  \Phi(X) \otimes \Phi(X) + D_6 D_5 L_f \Phi(X) \otimes \Phi(\tilde{X}) \\
& & \qquad + D_5 D_6 L_f \Phi(\tilde{X}) \otimes \Phi(X) + D_6 D_6 L_f \Phi(\tilde{X}) \otimes \Phi(\tilde{X}) \biggr] \\
& & + (\varepsilon - \tilde{\varepsilon})(1- \tilde{\varepsilon}-\varepsilon) \Ex_{\Q\otimes \P} 
    \biggl[D_5 D_5 L_f  \Phi(X) \otimes \Phi(X) + D_6 D_5 L_f \Phi(X) \otimes \Phi(\tilde{X}) \\
& & \qquad + D_5 D_6 L_f \Phi(\tilde{X}) \otimes \Phi(X) + D_6 D_6 L_f \Phi(\tilde{X}) \otimes \Phi(\tilde{X}) \biggr] \\
& & + (\varepsilon - \tilde{\varepsilon})(\tilde{\varepsilon}+\varepsilon) \Ex_{\Q\otimes \Q} 
\biggl[D_5 D_5 L_f  \Phi(X) \otimes \Phi(X) + D_6 D_5 L_f \Phi(X) \otimes \Phi(\tilde{X}) \\
& & \qquad + D_5 D_6 L_f \Phi(\tilde{X}) \otimes \Phi(X) + D_6 D_6 L_f \Phi(\tilde{X}) \otimes \Phi(\tilde{X}) \biggr] \,.
\end{eqnarray*}

Let $\varepsilon,\tilde{\varepsilon}\in[-1,+1]$. 
Then assumption {(\ref{loss-assump2})} and the bound $\|\Phi(\tilde{X}) \otimes \Phi(\tilde{X})\|_{{\mathcal L}(H, H)}\leq \inorm{k}^2$ 
yield the following inequality for the norm of this difference: 
$$  \left\|\frac{\partial G}{\partial H}(\varepsilon, \tilde{f}) - \frac{\partial G}{\partial H}(\tilde{\varepsilon}, \tilde{f})\right\|_{{\mathcal L}(H, H)} 
\leq 
4 c_{L, 2} \inorm{k}^2 
|\varepsilon - \tilde{\varepsilon}| \left(4 + 4 |\varepsilon| + 4 |\tilde{\varepsilon}|\right). $$
 Thus we have 
\begin{eqnarray*} 
& & \left\|\frac{\partial G}{\partial H}(\varepsilon, f) - \frac{\partial G}{\partial H}(\tilde{\varepsilon}, \tilde{f})\right\|_{{\mathcal L}(H, H)} \\
& \leq & 4 \inorm{k}^2 \omega (\inorm{k} \|f- \tilde{f}\|_H)_{r \inorm{k}} + 4 c_{L, 2} \inorm{k}^2
|\varepsilon - \tilde{\varepsilon}| \left(4 + 4 |\varepsilon| + 4 |\tilde{\varepsilon}|\right). 
\end{eqnarray*} 
Then the continuity of the partial derivative $\frac{\partial G}{\partial H}$ follows. 
This proves the continuous differentiability of $G$. 

Let $f\in H$. Consider the linear operator $\frac{\partial G}{\partial H}(0, f)$. It can be expressed as 
\begin{eqnarray}
& & \frac{\partial G}{\partial H}(0, f)  \label{dGdH} \\
&=& 2 \lambda id_H + \Ex_{\P \otimes \P}
\Bigl[D_5 D_5 L_f \Phi(X) \otimes \Phi(X) + D_6 D_5 L_f \Phi(X) \otimes \Phi(\tilde{X}) \nonumber \\
&& \qquad + D_5 D_6 L_f \Phi(\tilde{X}) \otimes \Phi(X) + D_6 D_6 L_f \Phi(\tilde{X}) \otimes \Phi(\tilde{X}) \Bigr].
\nonumber 
\end{eqnarray}
It is important to note that by Assumption \ref{assumption-loss1} $L$ is a twice continuously differentiable pairwise loss function which implies 
that 
$$
D_6 D_5 L_f = D_5 D_6 L_f \, , \qquad f\in H.
$$ 
Hence for any $g, \tilde{g} \in H$, the following holds 
\begin{eqnarray*}
& & \Bigl\langle \frac{\partial G}{\partial H}(0, f)(g), \tilde{g}\Bigr\rangle_H \\
& = & 2 \lambda \langle g, \tilde{g}\rangle_H + \Ex_{\P \otimes \P}
\biggl[D_5 D_5 L_f  \, g(X) \tilde{g}(X) \\
& & + D_6 D_5 L_f \, \Bigl(g(X) \tilde{g}(\tilde{X})+ g (\tilde{X}) g(\tilde{X})\Bigr) + D_6 D_6 L_f \, g(\tilde{X}) \tilde{g}(\tilde{X}) \biggr].
\end{eqnarray*}
So the linear operator $\frac{\partial G}{\partial H}(\varepsilon, \tilde{f})(0, f)$ is symmetric. Hence its spectrum lies in the closed interval $[a, b]$ where 
 $$ 
 a:=\inf_{\|g\|_H =1}\Bigl\langle \frac{\partial G}{\partial H}(0, f)(g), g\Bigr\rangle_H, \quad 
 b:=\sup_{\|g\|_H =1}\Bigl\langle \frac{\partial G}{\partial H}(0, f)(g), g\Bigr\rangle_H. 
 $$
Now, $L$ is  a \emph{convex} pairwise loss function due to Assumption \ref{assumption-loss2}. Therefore, we obtain, for any $g\in H$, 
$$ 
\Bigl \langle \frac{\partial G}{\partial H}(0, f)(g), g \Bigr \rangle_H \geq 2 \lambda \|g\|_H^2. $$ 
Hence, $a \geq 2 \lambda >0$. This shows that the operator $\frac{\partial G}{\partial H}(\varepsilon, \tilde{f})(0, f)$ is invertible. 
\qedr
\end{proofof}

We are now ready for the

\begin{proofof}{\textbf{Proof of Theorem \ref{thm.representer}}}
We will first prove part \emph{(i)} using Lemma \ref{Lemma-root}.
Let $\P\in\PM(\cXY)$.
Fix some $\xi\in\cX$ and define $g_\xi:=\Phi(\xi)=k(\cdot,\xi)\in H$.
By the reproducing property {(\ref{kernel-prop0})} of the kernel $k$, we have
\be \label{proof.representerthm.f0}
\langle \fPLs, g_\xi \rangle_H 
=  \langle \fPLs, \Phi(\xi) \rangle_H
= \fPLs(\xi).
\ee
Obviously, we also have
$g_\xi(x)=\Phi(\xi)(x)=k(x,\xi)$ and 
$g_\xi(\tix)=\Phi(\xi)(\tix)=k(\tix,\xi)$ for all $x,\tix\in\cX$.
Note that the partial derivatives of $L$ and of $\Ls$ with respect
of the last two arguments are identical, because $L$ and its shifted version $\Ls$ differ only by the term $L(\xyxy,0,0)$ which is independent of $f\in H$.
Therefore, Lemma \ref{Lemma-root} yields for the function $g_\xi\in H$
the equality
\begin{eqnarray*}
0 & = & \Ex_{\P^2} \Big[ \DfiveL{\fPLs} g_\xi(X) \\
    &   & ~~~~~~ + \DsixL{\fPLs} g_\xi(\tilde{X})\Big] 
                        + 2\lb \langle \fPLs, g_\xi \rangle_H \\
    &  = &  \Ex_{\P^2} \big[ h_{5,\P}(\XYXY) k(X,\xi)  + h_{6,\P}(\XYXY) k(\tilde{X},\xi) \big] 
                        + 2\lb \fPLs(\xi) ,
\end{eqnarray*}
where we used in the last equality the definition of $h_{5,\P}$ and $h_{6,\P}$
from {(\ref{thm.representer.f2})} and {(\ref{thm.representer.f3})}, respectively,
and {(\ref{proof.representerthm.f0})}.
From this we easily conclude that, for \emph{all} $\xi \in \cX$, 
\begin{eqnarray*}
  \fPLs(\xi) & = & -\frac{1}{2\lb}  \Ex_{\P^2} \big[h_{5,\P}(\XYXY) k(X,\xi) + 
                            h_{6,\P}(\XYXY) k(\tilde{X},\xi) \big] 
\end{eqnarray*}
which gives the assertion of part \emph{(i)}.

Let us now prove part \emph{(ii)}.
As $\Ls$ and $\lb\in(0,\infty)$ are fixed, we
will use the abbreviations $f_\P:=\fPLs$ and $f_\Q:=\fQLs$ in the proof.
The inequality is trivial, if $f_\P = f_{\Q}$. 
Hence let us assume that $f_{\P}\ne f_{\Q}$.
Recall the following well-known inequality from convex analysis. If $g:\R^2 \to \R$ is a convex and
total differentiable function with continuous second partial derivatives, then 
$$
   g(\tilde{t})-g(t)  \ge   \langle \nabla g(t), \tilde{t}-t \rangle_{\R^2} 
  \qquad \forall \, t,\tilde{t}\in\R^2,
$$
where $\nabla g(t)$ is the gradient of $g$ at $t$,
see \citet[Thm. 25.1, p.\,242]{Rockafellar1970} for a more general result using subgradients. 
To apply this result, we define, for any fixed $(\xyxy)\in(\cXY)^2$, the function
$g: \R^2\to\R$, where $g(t_1,t_2):=\Ls(\xyxy,t_1,t_2)$ and $t:=(t_1,t_2)\in\R^2$.
For any $\tilde{t}:=(\tilde{t}_1,\tilde{t}_2)\in\R^2$ we thus obtain
\begin{eqnarray}
  & & \Ls(\xyxy,\tilde{t}_1,\tilde{t}_2) -  \Ls(\xyxy,t_1, t_2)   \label{proof.representerthm.f1}\\
  & \ge & D_5 \Ls(\xyxy,t_1,t_2) (\tilde{t}_1-t_1) + D_6 \Ls(\xyxy,t_1,t_2) (\tilde{t}_2-t_2) \,.
  \nonumber 
\end{eqnarray}
If we specialize $(t_1,t_2):=(f_\P(x),f_\P(\tix))$ and 
$(\tilde{t}_1,\tilde{t}_2):=(f_\Q(x),f_\Q(\tix))$, we obtain from {(\ref{proof.representerthm.f1})} the inequality
\begin{eqnarray*}
  &  & \Ls(\xyxy,f_\Q(x),f_\Q(\tix)) -  \Ls(\xyxy,f_\P(x), f_\P(\tix)) \\
  & \ge  & D_5 \Ls(\xyxy,f_\P(x),f_\P(\tix)) (f_\Q(x)-f_\P(x)) \\
  &  &  + \,D_6 \, \Ls(\xyxy,f_\P(x),f_\P(\tix)) (f_\Q(\tix)-f_\P(\tix))    \\
  & =  & D_5 L(\xyxy,f_\P(x),f_\P(\tix)) (f_\Q(x)-f_\P(x)) \\
  &  &  + \,D_6 L(\xyxy,f_\P(x),f_\P(\tix)) (f_\Q(\tix)-f_\P(\tix)) ,   \\   
\end{eqnarray*}
where we used in the last step that $L$ and $\Ls$ differ only a term which does not dependent on the last two arguments.
By calculating the corresponding Bochner integral with respect to the product measure $\Q^2$, 
it follows from the reproducing property {(\ref{kernel-prop0})} of  $k$ that
\begin{eqnarray}
  &  & \RLs{\Q}{f_\Q} - \RLs{\Q}{f_\P} \label{proof.representerthm.f3}\\
  &  \ge &  \Big\langle f_\Q - f_\P ~,~ \Ex_{\Q^2} \Big[ \DfiveL{f_\P} \Phi(X) \nonumber \\
  &        & ~~~~~~~~~~~~~~~~~~~~~ + \DsixL{f_\P} \Phi(\tiX)\Big] \Big\rangle_H \nonumber \\
 &  = &  \Big\langle f_\Q - f_\P ~,~ \Ex_{\Q^2} \big[ h_{5,\P}(\XYXY) \Phi(X)  + 
              h_{6,\P}(\XYXY) \Phi(\tiX)\big] \Big\rangle_H \,,    \label{proof.representerthm.f3b}
\end{eqnarray}
where we used in the last step only the definition of $h_{5,\P}$ and 
$h_{6,\P}$ given in {(\ref{thm.representer.f2})} and {(\ref{thm.representer.f3})}, 
respectively.
Moreover, an easy calculation shows
\begin{eqnarray} \label{proof.representerthm.f4}
  2 \lb \langle f_\Q - f_\P, f_\P \rangle_H + \lb \hhnorm{f_\P-f_\Q}  
  =  \lb \hhnorm{f_\Q} - \lb \hhnorm{f_\P} \, .
\end{eqnarray}
We thus obtain
\begin{eqnarray}
  &  &  \Big\langle f_\Q - f_\P, \Ex_{\Q^2} \big[ h_{5,\P}(\XYXY) \Phi(X)  + h_{6,\P}(\XYXY) \Phi(\tiX)\big] + 2\lb f_\P \Big
  \rangle_H \nonumber \\
  &         & + \lb \hhnorm{f_\P-f_\Q}  \label{proof.representerthm.f5} \\
  & \stackrel{{\scriptsize{(\ref{proof.representerthm.f3b}),(\ref{proof.representerthm.f3})}}}{\le}  & \RLs{\Q}{f_\Q} - \RLs{\Q}{f_\P} 
                + 2\lb \langle f_\Q -f_\P, f_\P\rangle_H + \lb \hhnorm{f_\P-f_\Q} \nonumber \\
 & \stackrel{{\scriptsize{(\ref{proof.representerthm.f4})}}}{=} & \RLs{\Q}{f_\Q} - \RLs{\Q}{f_\P} 
                + \lb \hhnorm{f_\Q} - \lb \hhnorm{f_\P} \nonumber \\
 & = & \RLsreg{\Q}{f_\Q} - \RLsreg{\Q}{f_\P} ~~ \le 0, \label{proof.representerthm.f6}
\end{eqnarray}
where the term on the left hand side of {(\ref{proof.representerthm.f6})} is less than or equal to zero, because the regularized risk with respect to $\Q$ is minimized for $f_\Q$.
Recall that we have  
\be
  2 \lb f_\P = - \Ex_{\P^2} \Big[ h_{5,\P}(\XYXY) \Phi(X)  + h_{6,\P}(\XYXY) \Phi(\tiX)\Big] \label{proof.representerthm.f7}
\ee
due to {(\ref{thm.representer.f1})} in the first part of the representer theorem.
If we combine these two facts with the Cauchy-Schwarz inequality we obtain
from  {(\ref{proof.representerthm.f5})} that 
\begin{eqnarray*}
&     &    \lb \hhnorm{f_\P - f_\Q}  \\
& \le &  -  ~ \Big \langle f_\Q - f_\P ~,~ \Ex_{\Q^2} \big[ h_{5,\P}(\XYXY) \Phi(X) + h_{6,\P}(\XYXY) \Phi(\tiX)\big] + 2\lb f_\P \Big \rangle_H\\
& = &  \Big \langle f_\P - f_\Q ~,~ \Ex_{\Q^2} \big[ h_{5,\P}(\XYXY) \Phi(X) + h_{6,\P}(\XYXY) \Phi(\tiX)\big] + 2\lb f_\P \Big \rangle_H\\
& \stackrel{{\scriptsize{(\ref{proof.representerthm.f7})}}}{=}   &   \Big \langle f_\P - f_\Q ~,~ \Ex_{\Q^2} \big[ h_{5,\P}(\XYXY) \Phi(X)  + h_{6,\P}(\XYXY) \Phi(\tiX)\big]   \\
  &         & ~~~~~~~~~~~~~ - \Ex_{\P^2} \big[ h_{5,\P}(\XYXY) \Phi(X)  + h_{6,\P}(\XYXY) \Phi(\tiX)\big]\Big \rangle_H \\
  & \stackrel{C.-S.}{\le} & \hnorm{f_\P - f_\Q} \cdot \Big \| \Ex_{\Q^2} \big[ h_{5,\P}(\XYXY)\Phi(X)  + h_{6,\P}(\XYXY) \Phi(\tiX)\big] \\
  &        & ~~~~~~~~~~~~~~~~~~~ - \Ex_{\P^2} \big[ h_{5,\P}(\XYXY)\Phi(X)  + h_{6,\P}(\XYXY) \Phi(\tiX)\big]\Big \|_H ~. \\
\end{eqnarray*}
After multiplication with $1/\bigl( \lb \hnorm{f_\P - f_\Q} \bigr)$, which is allowed since $f_{\P} \ne f_{\Q}$, we immediately obtain the assertion.
\qedr
\end{proofof}

\begin{proofof}{\textbf{Proof of Theorem \ref{boundedmaxbias}}}
The proof only needs some elementary arguments. 
For brevity let us use the notation $\P_\e:=(1-\e)\P+\e\bar{\P}$,
$f_\P:=\fPLs$, $f_\Q:=\fQLs$, $f_{\P_\e}:=\fPepsLs$, and
$L_f:= L(X,Y,\tiX,\tiY,f(X),f(\tiX))$, $f\in H$.
Denote $f_0:=0\in H$.
Because $L\in[0,c]$ and $\lb \hhnorm{f_0}=0$, we immediately obtain, for all $\P\in\PM(\cZ)$,
$$
0 \le R^{reg}(\P) =\inf_{f\in H} \Big( \int L_{f} \,d\P^2 + \lb \hhnorm{f} \Bigr)
  \le \int L_{f_0} \,d\P^2 + \lb \hhnorm{f_0} \le c.
$$
Hence, there is no need to consider shifted loss functions.

Let us start with  part \emph{(i)}.
Because $f_\Q\in H$ and $L\in[0,c]$, we have $R^{reg}(\Q) = \int L_{f_\Q} \,d\Q^2 + \lb \hhnorm{f_\Q}$ and
$R^{reg}(\P) \le \int L_{f_\Q} \,d\P^2 + \lb \hhnorm{f_\Q}$.
Therefore,
\begin{eqnarray*}
 R^{reg}(\Q)-R^{reg}(\P) & \ge & \int L_{f_\Q} \,d\Q^2 + \lb \hhnorm{f_\Q} - \int L_{f_\Q} \,d\P^2 - \lb \hhnorm{f_\Q} \\
 & = & \int L_{f_\Q} \,d\Q^2 - \int L_{f_\Q} \,d\P^2 \\
 & \ge & -c \, d_{TV}(\Q^2,\P^2) 
  \ge  -2c \,  d_{TV}(\Q,\P),
\end{eqnarray*}
where we used in the last inequality that 
\begin{equation}\label{HoeffdingWolfowitzInequality}
 d_{TV}(\Q,\P) \le d_{TV}(\Q^2,\P^2) \le 2 d_{TV}(\Q,\P) , \quad \P,\Q\in\PM(\cXY),
\end{equation}
see \citet[p.709, (4.4) and (4.5)]{HoeffingWolfowitz1958}.
Analogously, from $f_\P\in H$ and $L(\xyxytt)\in[0,c]$, we conclude $R^{reg}(\Q) \le \int L_{f_\P} \,d\Q^2 + \lb \hhnorm{f_\P}$ and
$R^{reg}(\P) = \int L_{f_\P} \,d\P^2 + \lb \hhnorm{f_\P}$, which yields 
\begin{eqnarray*}
 R^{reg}(\Q)-R^{reg}(\P) & \le & \int L_{f_\P} \,d\Q^2 + \lb \hhnorm{f_\P} - \int L_{f_\P} \,d\P^2 - \lb \hhnorm{f_\P} \\
 & = & \int L_{f_\P} \,d\Q^2 - \int L_{f_\P} \,d\P^2 \\
 & \le & c \, d_{TV}(\Q^2,\P^2)  \stackrel{\scriptsize{(\ref{HoeffdingWolfowitzInequality})}}{\le} 2c \, d_{TV}(\Q,\P).
\end{eqnarray*}
If we combine both inequalities, we obtain the assertion from part \emph{(i)}.

To part \emph{(ii)}.
Because $f_{\P_\e}\in H$ and $L(\xyxytt)\in[0,c]$, we have $R^{reg}(\P_\e) = \int L_{f_{\P_\e}} \,d\P_\e^2 + \lb \hhnorm{f_{\P_\e}}$ and
$R^{reg}(\P) \le \int L_{f_{\P_\e}} \,d\P^2 + \lb \hhnorm{f_{\P_\e}}$.
Therefore,
\begin{eqnarray*}
 & & R^{reg}(\P_\e)-R^{reg}(\P) \\
 & \ge & \int L_{f_{\P_\e}} \,d\P_\e^2 - \int L_{f_{\P_\e}} \,d\P^2 \\
 & = & \int \int L_{f_{\P_\e}} \,d((1-\e)\P+\e\bar{\P}) \, d((1-\e)\P+\e\bar{\P}) - \int \int L_{f_{\P_\e}} \,d\P\,d\P\\
 & = & \left( (1-\e)^2-1\right) \int \int  L_{f_{\P_\e}} \,d\P\,d\P + \e(1-\e) \int \int  L_{f_{\P_\e}} \,d\P\,d\bar{\P}  \\
 &   & + \e(1-\e) \int \int  L_{f_{\P_\e}} \,d\bar{\P}\,d\P + \e^2 \int \int  L_{f_{\P_\e}} \,d\bar{\P}\,d\bar{\P} \\
 & = & \e \left( \int \int  L_{f_{\P_\e}} \,d\P\,d\bar{\P} + \int \int  L_{f_{\P_\e}} \,d\bar{\P}\,d\P 
        - 2 \int \int  L_{f_{\P_\e}} \,d\P\,d\P\right) \\
 &   & + \e^2 \left( \int \int  L_{f_{\P_\e}} \,d\P\,d\P + \int \int  L_{f_{\P_\e}} \,d\bar{\P}\,d\bar{\P} 
                     - \int \int  L_{f_{\P_\e}} \,d\bar{\P}\,d\P - \int \int  L_{f_{\P_\e}} \,d\P\,d\bar{\P}\right) \\
 & = & \e \left( \int \left[ \int  L_{f_{\P_\e}} \,d\P \right] \,d\bigl(\bar{\P}-\P\bigr)  + 
                 \int \left[ \int  L_{f_{\P_\e}} \,d\bigl( \bar{\P}-\P \bigr) \right] \,d\P \right) \\
 &   & + \e^2 \left( \int \left[ \int  L_{f_{\P_\e}} \,d\P\right] \,d\bigl(\P- \bar{\P}\bigr) 
                     + \int \left[\int L_{f_{\P_\e}} \,d\bar{\P}\right] \,d\bigl(\bar{\P}-\P\bigr)\right) \\ 
 & \stackrel{\scriptsize{(*)}}{\ge} & -2c\, d_{TV}(\bar{\P},\P) \e - 2c\, d_{TV}(\bar{\P},\P)\e^2\\
 & =   & -2c\, d_{TV}(\bar{\P},\P) \e (1+\e),
\end{eqnarray*}
where we used in {(*)} that $\int  L_{f_{\Q}} \,d\Q\le c$ for all $\Q\in\PM(\cXY)$.
Because $f_{\P}\in H$ and $L(\xyxytt)\in[0,c]$, we have $R^{reg}(\P_\e) \le \int L_{f_\P} \,d\P_\e^2 + \lb \hhnorm{f_\P}$ and
$R^{reg}(\P) = \int L_{f_\P} \,d\P^2 + \lb \hhnorm{f_\P}$. Hence, we obtain with the same argumentation as given above that
\begin{eqnarray*}
R^{reg}(\P_\e)-R^{reg}(\P) & \le & 2c\, d_{TV}(\bar{\P},\P) \e (1+\e).
\end{eqnarray*}
The combination of both inequalities yields the assertion.
\qedr
\end{proofof}

\begin{proofof}{\textbf{Proof of Theorem \ref{thm-boundedderivative}}}
The proof uses similar arguments than the proof of 
Theorem 15 in \citet{ChristmannSteinwart2007a}.

Partial derivatives of $L$ with respect to the
fifth or sixt argument are denoted by $D_5 L$ or $D_6 L$, respectively.
In the same manner we denote partial derivatives of $L$ of order two by
$D_i D_j L$, where $i,j\in\{5,6\}$.
Recall that due to {(\ref{frecheta1})}, $D_i\Ls=D_iL$ and $D_i D_j\Ls=D_i D_j L$ for 
$i,j\in\{5,6\}$.

Fix $\Q\in\PM(\cXY)$ and $\lb\in(0,\infty)$.
Define $\Pe:=(1-\e)\P+\e\Q$ and denote the product measure $\Pe\otimes\Pe$ by 
$\Pe^2$.

The function $G:\R \times H\to H$ defined by
\beqna
G(\e,f) & := &   2\lb f + \Ex_{\Pe^2} \Bigl[ D_5 L(X,Y,\tiX,\tiY,f(X),f(\tiX))\fm(X)   \\
&  &      \qquad \qquad \quad +     D_6 L(X,Y,\tiX,\tiY,f(X),f(\tiX))\fm(\tiX)\Bigr],
\eeqna
where $\e \in \R$ and $f \in H$, plays a key role in the proof.
Since $k$ is bounded by Assumption \ref{assumption-kernel1}, we have $\inorm{f}<\infty$ for all $f\in H$, see
{(\ref{kernel-prop1})}.
Additionally the partial derivatives $D_5 L$ and $D_6 L$ are continuous and uniformly
bounded by Assumption \ref{assumption-loss1}. Hence we obtain by using $\hnorm{ \int f\,d\P} \le \int \hnorm{f} \,d\P$
and {(\ref{kernel-prop3})}, that,
for all $\e\in\R$ and all $f\in H$,
\beqna
& & \hnorm{G(\e,f)} \\
& \le & 2 \lb \,\hnorm{f}  \\
& & + \int_{(\cXY)^2} \Bigl( \bigl( | D_5 L(\xyxy,f(x),f(\tix)| + |D_6 L(\xyxy,f(x),f(\tix))|\bigr) \cdot \\
& & ~~~~~~~~ \cdot \sup_{x\in \cX} \hnorm{\fm(x)} \Bigr)\,d\Pe^2(\xyxy) \\
& \stackrel{\footnotesize (\ref{kernel-prop2})}{\le} & 2 \lb \,\hnorm{f} + 2 c_{L,1} \cdot \inorm{k} < \infty\, .
\eeqna
Therefore, the map $G$ is well-defined and bounded with respect to the $H$-norm.
Due to {(\ref{kernel-prop1})} we have
$$
\inorm{G(\e,f)} \le 2 \bigl( \lb \,\hnorm{f} + c_{L,1} \cdot \inorm{k} \bigr) \, \inorm{k} ~ < ~ \infty\,.
$$
Note that for $\e\not\in [0,1]$ the $H$-valued Bochner integral is with respect to a signed measure.
Now for $\e\in[0,1]$ we obtain by using Lemma \ref{loss:diff-risk} that 
\begin{equation}\label{help-1}
G(\e,f) = \frac {\partial  (\cR_{\Ls,\Pe}(\cdot)+\lb\hhnorm{\cdot})}{\partial H}(f)\, .
\end{equation}
Since the map $f\mapsto \cR_{\Ls,\Pe}(f)+\lb\hhnorm{f}$
is strictly convex for all $\e \in[0,1]$ due to Lemma \ref{loss:convex-loss-lemma} and
is continuous,
equation (\ref{help-1}) shows that we have $G(\e,f)=0$
if and only if $f=f_{\Ls,\Pe,\lb}$ for such $\e$.
Our aim is to show the existence of a differentiable
function $\e\mapsto f_\e$ defined on a small interval $(-\d,\d)$ for some $\d>0$ 
that satisfies $G(\e,f_\e)=0$ for all
$\e\in (-\d,\d)$. Once we have shown the existence of this function we immediately obtain
\be \label{DSPQ}
S'_G(\P)(\Q) = \frac {\partial f_{\e}} {\partial \e} (0)\, .
\ee
For the existence of this map $\e\mapsto f_\e$ we have to check by the implicit function
theorem (\emph{cf.} Theorem \ref{implicit-thm})
that $G$ is \emph{continuously differentiable}
and that $\frac{\partial G}{\partial H}(0,f_{\P,\lb})$ is \emph{invertible}.
However, these properties of $G$ were shown in Theorem \ref{THM.PropertiesOfG}.
Hence we can apply Theorem \ref{implicit-thm} on implicit functions
to see that the map  $\e\mapsto f_\e$ is differentiable
on a small non-empty interval $(-\d,\d)$. 
Therefore, we obtain 
\begin{eqnarray*}
S'_G(\P)(\Q)
  & \stackrel{\footnotesize{{(\ref{DSPQ})}}}{=} &
   \frac {\partial f_{\e}} {\partial \e} (0) \\
  & \stackrel{\footnotesize{{(\ref{implicit-thm-f1})}}}{=} &
  -\left( \frac{\partial G}{\partial H}(0,\fPLs)\right)^{-1} \circ \frac{\partial G}{\partial \e}(0,\fPLs) \\  
  & \stackrel{\footnotesize{{(\ref{dGdeps}),~ (\ref{dGdH})}}}{=} & 
    -M(\P)^{-1} T(\Q;\P)\, ,
\end{eqnarray*}
which yields the assertion.
\qedr
\end{proofof}

\begin{proofof}{\textbf{Proof of Corollary \ref{cor-boundedIF}}}
The proof follows immediately by specifying $\Q$ to the Dirac-measure
$\delta_{(x_0,y_0)}$ in Theorem \ref{thm-boundedderivative}.
\qedr
\end{proofof}

\begin{proofof}{\textbf{Proof of Theorem \ref{thm.continuityofoperator}}}
We will first prove part \emph{(i)}. 
Let $\P\in\PM(\cXY)$ be fixed. Because 
$\Ls$ and $\lb$ are fixed, we will use again the shorter notations
$$
   h_{i,\P}(\XYXY)   :=  D_i\Ls\bigl(X,Y,\tiX,\tiY,\fPLs(X),\fPLs(\tiX)\bigr) , \quad i\in\{5,6\}, 
$$
in the proof, see {(\ref{thm.representer.f2})} and {(\ref{thm.representer.f3})}.

Let $\P_n\in\PM(\cXY)$, $n\in\N$, be a weakly convergence sequence with limit $\P$, i.e. $\P_n  \rightsquigarrow \P$. We know from {(\ref{WeakConvergenceAndDBL})} that $\P_n \rightsquigarrow \P$ is equivalent to $d_{\mathrm{BL}}(\P_n,\P)\to 0$, where $\dBL$ denotes the 
bounded Lipschitz metric, because $\cXY$ is separable by Assumption \ref{assumption-spaces1}.
Hence the metric space $(\cXY)^2$ is separable, too.
The separability guarantees that 
\be \label{proof.thm.continuityofoperator.f1}
  \Pn \rightsquigarrow \P \qquad  \Longleftrightarrow \qquad 
  \P_n^2 \rightsquigarrow \P^2 \qquad (n\to\infty), 
\ee
see \citet[Thm.\,2.8\,(ii), p.\,23]{Billingsley1999}.
Note that $\P_n^2 \rightsquigarrow \P^2$ guarantees by definition the convergence
$\int g\,d\P_n^2 \to \int g\,d\P^2$ for all continuous and bounded \emph{real-valued} functions $g:(\cXY)^2\to\R$. 
However, we will need a corresponding convergence result of Bochner integrals where the integrand is a special 
\emph{$H$-valued} function.
 
The second part of Theorem \ref{thm.representer} (representer theorem) yields
\begin{eqnarray}
 & & \hnorm{S(\P_n) - S(\P)}   :=  \hnorm{f_{\Ls,\P_n,\lb} - \fPLs}\label{proof.thm.continuityofoperator.f3} \\
 & \stackrel{\scriptsize{{(\ref{thm.representer.f4})}}}{\le} & \frac{1}{\lb} \Big\| 
                      \int  
                       \big[ 
                         h_{5,\P}(\XYXY) \Phi(X)  + h_{6,\P}(\XYXY) \Phi(\tiX)
                       \big]  \,d\P_n^2  \label{proof.thm.continuityofoperator.f4} \\
           & & ~~~- \int
                       \big[ 
                         h_{5,\P}(\XYXY) \Phi(X)  + h_{6,\P}(\XYXY) \Phi(\tiX)
                       \big]  \,d\P^2 \Big\|_H \, , \nonumber
\end{eqnarray}
where $\Phi(X)=k(\cdot,X)$ and $\Phi(\tilde{X})=k(\cdot,\tilde{X})$.
Because $k$ is continuous and bounded by Assumption
\ref{assumption-kernel1}, the canonical feature map $\Phi$ is continuous and bounded, too, see e.g. \citet[Lemma 4.23, Lemma 4.29]{SC2008}.
Furthermore, because the shifted loss function $\Ls$ is by Assumption 
\ref{assumption-loss1} twice continuously differentiable and the partial derivatives are uniformly bounded, it follows that, for every 
fixed $\P\in\PM(\cXY)$ and every fixed $\lb\in(0,\infty)$, the function
\begin{eqnarray} 
   & & \Psi_\P : \big((\cXY)^2, d_{(\cXY)^2}\big) \to (H,d_H), \nonumber \\
   & & \Psi_\P(\xyxy):= h_{5,\P}(\xyxy) \Phi(x) + h_{6,\P}(\xyxy) \Phi(\tilde{x})  \label{proof.thm.continuityofoperator.f5}
\end{eqnarray}
is \emph{continuous and bounded}, where $d_H(\cdot,\cdot):=\hnorm{\cdot-\cdot}$.
We mention that the $H$-valued function $\Psi_\P$ does not depend on $\P_n$.
Because $\Psi_\P$ is continuous and bounded, we obtain from \citet[p. III.40]{Bourbaki2004} 
the following convergence result for Bochner integrals:
\be \label{proof.thm.continuityofoperator.f6}
\P_n^2 \rightsquigarrow \P^2
\qquad \Longrightarrow \qquad 
\lim_{{n\to\infty} }\int \Psi_\P \,d\P_n^2  = \int \Psi_\P \,d\P^2 \,,
\ee
see also \citet[Thm. A.1, p. 1000]{HableChristmann2011}.
Combining {(\ref{proof.thm.continuityofoperator.f1})}--{(\ref{proof.thm.continuityofoperator.f6})},
we obtain that $\P_n \rightsquigarrow \P$, which is equivalent to 
$\dBL(\P_n,\P)\to 0$ by {(\ref{WeakConvergenceAndDBL})},  implies $\hnorm{S(\P_n)-S(\P)} \to 0$, which is the 
assertion of part \emph{(i)}.

The proof of the second part follows immediately from part \emph{(i)} 
and the fact that the inclusion $\mathrm{id}: H\to \mathcal{C}_b(\cX)$ is 
continuous and bounded, see e.g. \citet[Lemma 4.28]{SC2008}.
\qedr
\end{proofof}

\begin{proofof}{\textbf{Proof of Corollary \ref{cor.continuityofestimator}}}
Let $(D_{n,m})_{m\in\N}$ be a sequence  in $(\cXY)^n$ which converges to some
$D_{n,0}\in(\cXY)^n$, if $m\to \infty$. Then, the corresponding sequence of empirical measures
weakly converges, i.e. 
$\D_{n,m}  \rightsquigarrow \D_{n,0}$, if $m\to\infty$. Therefore, the assertion follows from
Theorem \ref{thm.continuityofoperator} and 
$f_{\Ls,\D_n,\lb} = S_n(D_n)$.
\qedr 
\end{proofof}

\begin{proofof}{\textbf{Proof of Theorem \ref{thm.qualitativerobust}}}
Fix $\lb\in(0,\infty)$. 
We will first prove part \emph{(i)}. For any $D_n\in (\cXY)^n$ denote its empirical measure by
$\D_n:=\frac{1}{n}\sum_{i=1}^n \delta_{(x_i,y_i)}$.
According to Corollary \ref{cor.continuityofestimator}, the functions 
\be \nonumber
S_n : \bigl((\cXY)^n,d_{(\cXY)^n}\bigr) \to (H,d_H), \quad S_n(D_n)=f_{\Ls,\D_n,\lb} 
\ee
are continuous and therefore measurable with respect to the corresponding Borel-$\sigma$-algebras for every $n\in\N$. The mapping
\be \label{def.S}
S : \big(\PM(\cXY),\dBL\big)  \to (H,d_H), \quad S(\P)=\fPLs ,
\ee
is a continuous operator due to Theorem \ref{thm.continuityofoperator}.
Furthermore, 
$$
S_n(D_n) = S(\D_n) \qquad \forall~ D_n\in(\cXY)^n~~\forall~n\in\N.
$$
Because $\cX$ is a separable metric space and the kernel $k$ is continuous, the RKHS $H$ of $k$ is
separable, see e.g. \citet[Lemma 4.33, p. 130]{SC2008}. Hence $(H,d_H)$ is a complete and separable metric space. 

Therefore, the sequence of \RPL estimators 
$(\fDDnLs)_{n\in\N}$, where $\DD_n:=\frac{1}{n}\sum_{i=1}^n \delta_{(X_i,Y_i)}$,
is qualitatively robust for all Borel probability measures $\P\in\PM(\cXY)$ according to \citet[Thm. 2]{Cuevas1988}, which states: If $(S_n)_{n\in\N}$ is any sequence of estimators which can be represented via a \emph{continuous} operator $S$, which maps each probability measure $\P$ to a value
in a \emph{complete and separable metric space} and satisfies (in our notation) $S_n(D_n)=S(\D_n)$, is qualitatively robust for all $\P$. Hence the assertion of part \emph{(i)} is shown. 

Let us now prove part \emph{(ii)}. 
It follows from the first part of Theorem \ref{thm.continuityofoperator}, that the operator 
$S$ defined in {(\ref{def.S})} is continuous for all $\P\in\PM(\cXY)$.
Hence all conditions of Assumption 16.3 in \citet{ChristmannSalibianBarreraVanAelst2013}
are satisfied, because
$\mathscr{Z}:=\cXY$ is a \emph{compact} metric space by assumption of 
Theorem \ref{thm.qualitativerobust}(ii) and 
$\mathscr{W}:=H$  is a complete and separable metric space due to the continuity of $k$ 
by Assumption \ref{assumption-kernel1},  e.g. \citet[Lemma 4.33, p. 130]{SC2008}.
Hence, Corollary 16.1 by  \citet{ChristmannSalibianBarreraVanAelst2013} is applicable and immediately yields the assertion.
We like to note that the compactness of the metric space $\mathscr{Z}$ was used in the proof of the above mentioned
Corollary 16.1 to show that the continuous operator $S$ is even uniformly continuous for all $\P\in\PM(\mathscr{Z})$. 
\qedr
\end{proofof}



\begin{thebibliography}{}

\bibitem[Agarwal and Niyogi(2009)Agarwal and Niyogi]{AgarwalNiyogi2009}
Agarwal, S. and Niyogi, P. (2009).
\newblock Generalization bounds for ranking algorithms via algorithmic
  stability.
\newblock {\em J. Mach. Learn. Res.}, {\bf 10}, 441--474.

\bibitem[Akerkar(1999)Akerkar]{Akerkar1999}
Akerkar, R. (1999).
\newblock {\em Nonlinear Functional Analysis\/}.
\newblock Narosa Publishing House, New Dehli.

\bibitem[Billingsley(1999)Billingsley]{Billingsley1999}
Billingsley, P. (1999).
\newblock {\em Convergence of {P}robability {M}easures\/}.
\newblock John Wiley {\&} Sons, New York, 2nd edition.

\bibitem[Bourbaki(2004)Bourbaki]{Bourbaki2004}
Bourbaki, N. (2004).
\newblock {\em Integration I. (Translated from the 1959, 1965, and 1967 French
  originals by Sterling K. Berberian. Chapters 1--6)\/}.
\newblock Springer, Berlin.

\bibitem[Cao {\em et~al.}(2015)Cao, Guo, and Ying]{CaoGuoYing2015}
Cao, Q., Guo, Z.~C., and Ying, T. (2015).
\newblock Generalization bounds for metric and similarity learning.
\newblock {\em to appear in: Machine Learning, Online First. DOI
  10.1007/s10994-015-5499-7\/}.

\bibitem[Castaing and Valadier(1977)Castaing and Valadier]{CaVa77}
Castaing, C. and Valadier, M. (1977).
\newblock {\em Convex Analysis and Measurable Multifunctions\/}.
\newblock Springer, Berlin.

\bibitem[Christmann and Steinwart(2004)Christmann and
  Steinwart]{ChristmannSteinwart2004a}
Christmann, A. and Steinwart, I. (2004).
\newblock On robust properties of convex risk minimization methods for pattern
  recognition.
\newblock {\em J. Mach. Learn. Res.}, {\bf 5}, 1007--1034.

\bibitem[Christmann and Steinwart(2007)Christmann and
  Steinwart]{ChristmannSteinwart2007a}
Christmann, A. and Steinwart, I. (2007).
\newblock Consistency and robustness of kernel based regression.
\newblock {\em Bernoulli\/}, {\bf 13}, 799--819.

\bibitem[Christmann {\em et~al.}(2009)Christmann, {V}an Messem, and
  Steinwart]{ChristmannVanMessemSteinwart2009}
Christmann, A., {V}an Messem, A., and Steinwart, I. (2009).
\newblock On consistency and robustness properties of support vector machines
  for heavy-tailed distributions.
\newblock {\em Statistics and Its Interface\/}, {\bf 2}, 311--327.

\bibitem[Christmann {\em et~al.}(2013)Christmann, Salib\'{i}an-Barrera, and
  Aelst]{ChristmannSalibianBarreraVanAelst2013}
Christmann, A., Salib\'{i}an-Barrera, M., and Aelst, S.~V. (2013).
\newblock Qualitative robustness of bootstrap approximations for kernel based
  methods.
\newblock In C.~Becker, R.~Fried, and S.~Kuhnt, editors, {\em Robustness and
  Complex Data Structures\/}, pages 277--293. Springer, Heidelberg.

\bibitem[Clemencon {\em et~al.}(2008)Clemencon, Lugosi, and
  Vayatis]{Clemencon2008}
Clemencon, S., Lugosi, G., and Vayatis, N. (2008).
\newblock Ranking and empirical minimization of {U}-statistics.
\newblock {\em Ann. Statist.}, {\bf 36}, 844--874.

\bibitem[Cucker and Zhou(2007)Cucker and Zhou]{CuckerZhou2007}
Cucker, F. and Zhou, D.~X. (2007).
\newblock {\em Learning Theory: An Approximation Theory Viewpoint\/}.
\newblock Cambridge University Press, Cambridge.

\bibitem[Cuevas(1988)Cuevas]{Cuevas1988}
Cuevas, A. (1988).
\newblock Qualitative robustness in abstract inference.
\newblock {\em J. Statist. Plann. Inference\/}, {\bf 18}, 277--289.

\bibitem[Denkowski {\em et~al.}(2003)Denkowski, Mig{\'o}rski, and
  Papageorgiou]{DenkowskiEtAl2003}
Denkowski, Z., Mig{\'o}rski, S., and Papageorgiou, N. (2003).
\newblock {\em An {I}ntroduction to {N}onlinear {A}nalysis: {T}heory\/}.
\newblock Kluwer Academic Publishers, Boston.

\bibitem[Dudley(2002)Dudley]{Dudley2002}
Dudley, R.~M. (2002).
\newblock {\em {R}eal {A}nalysis and {P}robability\/}.
\newblock Cambridge University Press, Cambridge.

\bibitem[Efron(1979)Efron]{Efron1979}
Efron, B. (1979).
\newblock Bootstrap methods: Another look at the jackknife.
\newblock {\em Ann. Statist.}, {\bf 7}, 1--26.

\bibitem[Efron(1982)Efron]{Efron1982}
Efron, B. (1982).
\newblock {\em The Jackknife, the Bootstrap, and Other Resampling Plans\/},
  volume~38.
\newblock CBMS Monograph, Society for Industrial and Applied Mathematics,
  Philadelphia.

\bibitem[Ekeland and T\'{e}mam(1999)Ekeland and T\'{e}mam]{EkelandTemam1999}
Ekeland, I. and T\'{e}mam, R. (1999).
\newblock {\em Convex Analysis and Variational Problems\/}.
\newblock SIAM, Philadelphia.

\bibitem[Ekeland and Turnbull(1983)Ekeland and Turnbull]{EkelandTurnbull1983}
Ekeland, I. and Turnbull, T. (1983).
\newblock {\em Infinite-Dimensional Optimization and Convexity\/}.
\newblock University of Chicago Press, Chicago.

\bibitem[Fan {\em et~al.}(2014)Fan, Hu, Wu, and Zhou]{FanHuWuZhou2013}
Fan, J., Hu, T., Wu, Q., and Zhou, D.~X. (2014).
\newblock Consistency analysis of an empirical minimum error entropy algorithm.
\newblock {\em to appear in: Appl. Comput. Harmonic Anal., Online first:
  \url{doi:10.1016/j.acha.2014.12.005}\/}.

\bibitem[Feng {\em et~al.}(2015)Feng, Huang, Shi, Yang, and
  Suykens]{FengHuangShiYangSuykens2015}
Feng, Y., Huang, X., Shi, L., Yang, Y., and Suykens, J. (2015).
\newblock Learning with the maximum correntropy criterion induced losses for
  regression.
\newblock {\em J. Mach. Learn. Res.}, {\bf 16}, 993--1034.

\bibitem[Hable and Christmann(2011)Hable and Christmann]{HableChristmann2011}
Hable, R. and Christmann, A. (2011).
\newblock Qualitative robustness of support vector machines.
\newblock {\em Journal of Multivariate Analysis\/}, {\bf 102}, 993--1007.

\bibitem[Hable and Christmann(2013)Hable and Christmann]{HableChristmann2013a}
Hable, R. and Christmann, A. (2013).
\newblock Robustness versus consistency in ill-posed classification and
  regression problems.
\newblock In A.~Giusti, G.~Ritter, and M.~Vichi, editors, {\em Classification
  and Data Mining\/}, pages 27--35. Springer, Berlin.

\bibitem[Hampel(1968)Hampel]{Hampel1968}
Hampel, F.~R. (1968).
\newblock Contributions to the theory of robust estimation.
\newblock Unpublished Ph.D. thesis, Department of Statistics, University of
  California, Berkeley.

\bibitem[Hampel(1971)Hampel]{Hampel1971}
Hampel, F.~R. (1971).
\newblock A general qualitative definition of robustness.
\newblock {\em Ann. Math. Statist.}, {\bf 42}, 1887--1896.

\bibitem[Hampel(1974)Hampel]{Hampel1974}
Hampel, F.~R. (1974).
\newblock The influence curve and its role in robust estimation.
\newblock {\em J. Amer. Statist. Assoc.}, {\bf 69}, 383--393.

\bibitem[Hampel {\em et~al.}(1986)Hampel, Ronchetti, Rousseeuw, and
  Stahel]{HampelRonchettiRousseeuwStahel1986}
Hampel, F.~R., Ronchetti, E.~M., Rousseeuw, P.~J., and Stahel, W.~A. (1986).
\newblock {\em Robust statistics: The Approach Based on Influence Functions\/}.
\newblock John Wiley {\&} Sons, New York.

\bibitem[Hoeffding and Wolfowitz(1958)Hoeffding and
  Wolfowitz]{HoeffingWolfowitz1958}
Hoeffding, W. and Wolfowitz, J. (1958).
\newblock Distinguishability of sets of distributions. the case of independent
  and identically distributed chance variables.
\newblock {\em Ann. Math. Statist.}, {\bf 29}, 700--718.

\bibitem[Hu {\em et~al.}(2013)Hu, Fan, Wu, and Zhou]{HuFanWunZhou2013}
Hu, T., Fan, J., Wu, Q., and Zhou, D.~X. (2013).
\newblock Learning theory approach to minimum error entropy criterion.
\newblock {\em J. Mach. Learn. Res.}, {\bf 14}, 377--397.

\bibitem[Hu {\em et~al.}(2015)Hu, Fan, Wu, and Zhou]{HuFanWuZhou2015}
Hu, T., Fan, J., Wu, Q., and Zhou, D.~X. (2015).
\newblock Regularization schemes for minimum error entropy principle.
\newblock {\em Anal. Appl.}, {\bf 13}, 437--455.

\bibitem[Huber(1967)Huber]{Huber1967}
Huber, P.~J. (1967).
\newblock The behavior of maximum likelihood estimates under nonstandard
  conditions.
\newblock {\em Proc. 5th Berkeley Symp.}, {\bf 1}, 221--233.

\bibitem[Huber(1981)Huber]{Huber1981}
Huber, P.~J. (1981).
\newblock {\em Robust Statistics\/}.
\newblock John Wiley {\&} Sons, New York.

\bibitem[Koroljuk and Borovskich(1994)Koroljuk and
  Borovskich]{KoroljukBorovskich1994}
Koroljuk, V. and Borovskich, Y. (1994).
\newblock {\em Theory of $U$-Statistics\/}.
\newblock Springer, Dordrecht.

\bibitem[Mukherjee and Zhou(2006)Mukherjee and Zhou]{MukherjeeZhou2006}
Mukherjee, S. and Zhou, D.~X. (2006).
\newblock Learning coordinate covariances via gradients.
\newblock {\em J. Mach. Learn. Res.}, {\bf 7}, 519--549.

\bibitem[Poggio and Girosi(1998)Poggio and Girosi]{PoggioGirosi1998}
Poggio, T. and Girosi, F. (1998).
\newblock A sparse representation for function approximation.
\newblock {\em Neural Comput.}, {\bf 10}, 1445--1454.

\bibitem[Principe(2010)Principe]{Principe2010}
Principe, J. (2010).
\newblock {\em Information Theoretic Learning: Renyi's Entropy and Kernel
  Perspectives\/}.
\newblock Springer, New York.

\bibitem[Rio(2013)Rio]{Rio2013}
Rio, E. (2013).
\newblock On {McDiarmid's} concentration inequality.
\newblock {\em Electron. Commun. Probab.}, {\bf 44}, 1--011.

\bibitem[Rockafellar(1970)Rockafellar]{Rockafellar1970}
Rockafellar, R.~T. (1970).
\newblock {\em Convex Analysis\/}.
\newblock Princeton University Press, Princeton, NJ.

\bibitem[Sch{\"o}lkopf and Smola(2002)Sch{\"o}lkopf and Smola]{ScSm2002}
Sch{\"o}lkopf, B. and Smola, A.~J. (2002).
\newblock {\em Learning with Kernels. Support Vector Machines, Regularization,
  Optimization, and Beyond\/}.
\newblock MIT Press, Cambridge, MA.

\bibitem[Serfling(1980)Serfling]{Serfling1980}
Serfling, R.~J. (1980).
\newblock {\em Approximation Theorems of Mathematical Statistics\/}.
\newblock John Wiley {\&} Sons, New York.

\bibitem[Steinwart and Christmann(2008)Steinwart and Christmann]{SC2008}
Steinwart, I. and Christmann, A. (2008).
\newblock {\em Support Vector Machines\/}.
\newblock Springer, New York.

\bibitem[Tukey(1977)Tukey]{Tukey1977}
Tukey, J.~W. (1977).
\newblock {\em Exploratory Data Analysis. [preliminary edition 1970-1971]\/}.
\newblock Addison-Wesley, Reading, MA.

\bibitem[Vapnik(1995)Vapnik]{Vapnik1995}
Vapnik, V.~N. (1995).
\newblock {\em The Nature of Statistical Learning Theory\/}.
\newblock Springer, New York.

\bibitem[Vapnik(1998)Vapnik]{Vapnik1998}
Vapnik, V.~N. (1998).
\newblock {\em Statistical Learning Theory\/}.
\newblock John Wiley {\&} Sons, New York.

\bibitem[Wahba(1999)Wahba]{Wahba1999}
Wahba, G. (1999).
\newblock Support vector machines, reproducing kernel {H}ilbert spaces and the
  randomized {GACV}.
\newblock In B.~Sch{\"o}lkopf, C.~J.~C. Burges, and A.~Smola, editors, {\em
  Advances in Kernel Methods--Support Vector Learning\/}, pages 69--88. MIT
  Press, Cambridge, MA.

\bibitem[Wendland(1995)Wendland]{Wendland1995}
Wendland, H. (1995).
\newblock Piecewise polynomial, positive definite and compactly supported
  radial basis functions of minimal degree.
\newblock {\em Adv. Comput. Math.}, {\bf 4}, 389--396.

\bibitem[Wu(1995)Wu]{Wu1995}
Wu, Z. (1995).
\newblock Compactly supported positive definite radial functions.
\newblock {\em Adv. Comput. Math.}, {\bf 4}, 283--292.

\bibitem[Xing {\em et~al.}(2002)Xing, Ng, Jordan, and
  Russell]{XinNgJordanRussell2002}
Xing, E., Ng, A., Jordan, M., and Russell, S. (2002).
\newblock Distance metric learning, with application to clustering with
  side-information.
\newblock {\em Advances in Neural Information Processing Systems\/}, {\bf 15},
  505--512.

\bibitem[Ying and Zhou(2015)Ying and Zhou]{YingZhou2015}
Ying, Y. and Zhou, D.~X. (2015).
\newblock Unregularized online learning algorithms with general loss functions.
\newblock {\em To appear in: Appl. Comput. Harmonic Anal., Online first:
  \url{doi:10.1016/j.acha.2015.08.007}\/}.

\end{thebibliography}

\end{document}